\documentclass[graybox]{svmult}
\usepackage{amsfonts}
\usepackage{amsmath}
\usepackage[english]{babel}
\ifx\pdftexversion\undefined
 \usepackage[dvipsone]{epsfig}
 \RequirePackage{graphics,epsf}
 
\else
 \usepackage{thumbpdf}
 \usepackage[pdftex]{graphics}
 \usepackage[pdftex,colorlinks=true,pdfstartview=FitV,linkcolor=black,citecolor=black,urlcolor=black,bookmarks=true]{hyperref}
\hypersetup{
       pdfauthor = {Ralf Deiterding},
       pdftitle = {},
       pdfsubject = {},
       pdfcreator = {LaTeX with hyperref package},
       pdfproducer = {dvips + ps2pdf}}
 
 \usepackage{graphicx}
\fi

\usepackage{mathptmx}       % selects Times Roman as basic font
\usepackage{helvet}         % selects Helvetica as sans-serif font
\usepackage{courier}        % selects Courier as typewriter font
\usepackage{type1cm}        % activate if the above 3 fonts are
\usepackage{makeidx}         % allows index generation
\usepackage{graphicx}        % standard LaTeX graphics tool
\usepackage{multicol}        % used for the two-column index
\usepackage[bottom]{footmisc}% places footnotes at page bottom

\allowdisplaybreaks

\begin{document}

\title*{Robust split-step Fourier methods for simulating the propagation of ultra-short pulses in single- and two-mode optical communication fibers}

\titlerunning{Robust split-step Fourier methods for propagation of ultra-short pulses in optical fibers} 

\author{Ralf Deiterding and Stephen W. Poole}
\institute{Ralf Deiterding \at German Aerospace Center (DLR), Institute of Aerodynamics and Flow Technology, 
Bunsenstr. 10, 37073 G\"ottingen, Germany, \email{ralf.deiterding@dlr.de}
\and Stephen W. Poole \at Oak Ridge National Laboratory, Computer Science and Mathematics Division, 
P.O. Box 2008 MS6164, Oak Ridge, TN 37831, USA, \email{spoole@ornl.gov}}
\maketitle

\abstract*{Extensions of the split-step Fourier method (SSFM) for Schr\"odinger-type pulse propagation equations for simulating
femto-second pulses in single- and two-mode optical communication fibers are developed and tested for Gaussian pulses. The core
idea of the proposed numerical methods is to adopt an operator splitting approach, in which the nonlinear sub-operator, consisting of 
Kerr nonlinearity, the self-steepening and stimulated Raman scattering terms, is reformulated using Madelung transformation into
a quasilinear first-order system of signal intensity and phase. A second-order accurate upwind numerical method is derived rigorously
for the resulting system in the single-mode case; a straightforward extension of this method is used to approximate the four-dimensional
system resulting from the nonlinearities of the chosen two-mode model. Benchmark SSFM computations of prototypical ultra-fast communication pulses
in idealized single- and two-mode fibers with homogeneous and alternating dispersion parameters and also high nonlinearity demonstrate
the reliable convergence behavior and robustness of the proposed approach.}

\abstract{Extensions of the split-step Fourier method (SSFM) for Schr\"odinger-type pulse propagation equations for simulating
femto-second pulses in single- and two-mode optical communication fibers are developed and tested for Gaussian pulses. The core
idea of the proposed numerical methods is to adopt an operator splitting approach, in which the nonlinear sub-operator, consisting of 
Kerr nonlinearity, the self-steepening and stimulated Raman scattering terms, is reformulated using Madelung transformation into
a quasilinear first-order system of signal intensity and phase. A second-order accurate upwind numerical method is derived rigorously
for the resulting system in the single-mode case; a straightforward extension of this method is used to approximate the four-dimensional
system resulting from the nonlinearities of the chosen two-mode model. Benchmark SSFM computations of prototypical ultra-fast communication pulses
in idealized single- and two-mode fibers with homogeneous and alternating dispersion parameters and also high nonlinearity demonstrate
the reliable convergence behavior and robustness of the proposed approach.}

\section{Introduction}
As computational capabilities are continuously rising, so is the demand for enhanced networking speed. One possible approach for increasing
data throughput is the design of networks with transmission speeds well in the Tb/s range. 
While the maximal single channel  communication speed in demonstrated wavelength division multiplexing systems is generally below 
$100\,\mathrm{Gb/s}$, cf. \cite{Gnauck-08}, we are in here concerned with the modeling of single- and two-mode mode optical fibers that are
suitable in particular for long-distance data transmission. 

At present, computational models for investigating the propagation of light pulses in fibers have been developed primarily for pulses with a temporal half 
width well in the pico-second regime. Pulses with half widths $T_0\gg 1\,\mathrm{ps}$ are sufficient for representing even on-off-key modulated bit streams with 
up to $100\,\mathrm{Gb/s}$ frequency. However, bit streams in the Tb/s regime can only be represented with {\em ultra-fast} pulses 
satisfying $T_0<100\,\mathrm{fs}$. Yet, in the ultra-fast pulse regime nonlinear pulse self-steepening and nonlinear stimulated Raman scattering
are not negligible anymore and an extended version of the Schr\"odinger-type pulse propagation equation has to be considered. 

Numerical solutions of the Schr\"odinger-type pulse propagation equation are primarily obtained with split-step Fourier schemes 
that perform spatial propagation steps considering firstly only the linearities in the equation by
discrete Fourier transformation and then secondly only the nonlinear terms. While the construction of such {\em split-step Fourier methods} 
(SSFM) is very well established, cf. \cite{Agrawal-07,Hohage-02}, the topic of how to incorporate both self-steepening and Raman scattering 
reliably into the SSFM has received little attention. Here, we will describe a new class of extended SSFM that properly consider
the hyperbolic nature of the nonlinear sub-operator for single- and coupled two-mode optical communication fibers.

The paper is organized as follows: In Sect.~\ref{sec:singlemodeeq}, we recall the governing equations of pulse propagation in single-mode 
fibers. Section~\ref{sec:smssfm} first discusses the construction principles of split-step Fourier methods and then proceeds by 
describing our new type of single-mode SSFM for ultra-fast pulses as first- and second-order accurate numerical schemes, cf. \cite{Deiterding-etal-13}. 
An ultra-fast Gaussian pulse benchmark confirming robust second-order accuracy of the overall SSFM and demonstrating its application for simulating pulse
propagation through an idealized dispersion-managed single-mode communication line are given. In Sect.~\ref{sec:twomode}, we 
describe an extended two-mode model for considering the simultaneous and fully coupled propagation of two ultra-fast pulses in a single fiber cable. 
The subsequent Sect.~\ref{sec:smsstm} presents a fractional step approach for effectively extending the derived single-mode nonlinear
sub-operator to the corresponding system in the two-mode case. A two-mode benchmark of two interacting ultra-fast Gaussian communication pulses confirms
the reliability of the method and its straightforward applicability in the dispersion-managed case is also shown. The conclusions are
given in Sect.~\ref{sec:conclusions}. 

\section{Governing equation for ultra-fast pulses in a single-mode fiber}\label{sec:singlemodeeq}
The most general equation representing single-mode pulse propagation in a one-dimensional optical fiber reads
\begin{equation}\label{eq:pulseRaman} \frac{\partial A}{\partial z}+\frac{\alpha}{2}A+\left(\sum_{k\ge 1}\beta_k \frac{i^{k-1}}{k!} 
\frac{\partial^k}{\partial t^k}\right)A = i \gamma\left(1+\frac{i}{\omega_0}\frac{\partial}{\partial t}\right)\times
\left[A\int^\infty_{-\infty} R(t')|A(t-t')|^2 dt'\right].  
\end{equation}
Equation~(\ref{eq:pulseRaman}) is derived from the electric field of the Maxwell equations, cf. \cite{Agrawal-07}, and describes the evolution of
the slowly varying field envelope $A(z,t)$ of the complex-valued signal over the propagation distance $z$ and time $t$. The coefficients $\beta_k$ 
model signal dispersion. Since the refractive index $n$ of the fiber material is dependent on the light's circular frequency $\omega$, 
different spectral components associated to a pulse travel at slightly different velocities, given by $c/n(\omega)$, with $c$ denoting the 
speed of light in vacuum. This effect is mathematically modeled by expressing the mode propagation constant $\beta$ in a Taylor series about the central 
frequency $\omega_0=2\pi c/\lambda_0$ as
\begin{equation}
\beta(\omega)=n(\omega)\frac{\omega}{c}=\sum_{k\ge 0}\frac{1}{k!} \beta_k (\omega-\omega_0)^k.
\end{equation}
Here, the wavelength of the injected laser light is denoted by $\lambda_0$ and the parameters $\alpha$ and $\gamma$ model linear signal loss and 
fiber nonlinearity, respectively. The function $R(t)$ represents intrapulse Raman scattering, a nonlinear effect transferring energy from higher 
to lower light frequencies. Using $R(t)=(1-f_R)\delta(t)+f_R h_R(t)$ with $f_R=0.18$ \cite{Blow-Wood-89} as Raman response function, 
applying a Taylor series expansion and neglecting higher order terms, Eq.~(\ref{eq:pulseRaman}) eventually becomes 
\begin{equation}\label{eq:pulsefast} 
\frac{\partial A}{\partial z}+\frac{\alpha}{2}A+\beta_1\frac{\partial A}{\partial t}+
i\frac{\beta_2}{2}\frac{\partial^2A}{\partial t^2}-\frac{\beta_3}{6}\frac{\partial^3A}{\partial t^3}=
i\gamma\left(A|A|^2 + \frac{i}{\omega_0}\frac{\partial}{\partial t}\left(A|A|^2\right)
-T_R A \frac{\partial|A|^2}{\partial t} \right).
\end{equation}
In general, Eq.~(\ref{eq:pulse}) is widely accepted as a valid model for modeling the propagation of pulses with a half 
width $T_0>10\,\mathrm{fs}$ \cite{Agrawal-07}. For $\lambda_0=1550\,\mathrm{nm}$, a typical value for the Raman response parameter
is $T_R=3\,\mathrm{fs}$. The first nonlinear term on the right-hand side of Eq.~(\ref{eq:pulsefast}) is called the {\em Kerr} nonlinearity
and the second represents nonlinear pulse self-steepening. 

Introducing the signal group velocity $v$ with $\beta_1=1/v$ and using the transformation $T\equiv t-z/v$ 
into {\em retarded} time $T$, Eq.~(\ref{eq:pulsefast}) is transformed into the frame of reference of the pulse to read
\begin{equation}\label{eq:pulse} 
\frac{\partial A}{\partial z}+\frac{\alpha}{2}A+
i\frac{\beta_2}{2}\frac{\partial^2 A}{\partial T^2}-\frac{\beta_3}{6}\frac{\partial^3 A}{\partial T^3}=
i\gamma\left(A|A|^2 + i S \frac{\partial}{\partial T}\left(A|A|^2\right)-T_R A \frac{\partial|A|^2}{\partial T} \right), 
\end{equation}
where we have also introduced $S=\omega_0^{-1}$. For $T_0\gg 1\,\mathrm{ps}$, the last two terms can be neglected and Eq.~(\ref{eq:pulse}) reduces to
\begin{equation}\label{eq:pulseslow} 
\frac{\partial A}{\partial z}+\frac{\alpha}{2}A+
i\frac{\beta_2}{2}\frac{\partial^2A}{\partial T^2}-\frac{\beta_3}{6}\frac{\partial^3A}{\partial T^3}=
i\gamma A|A|^2,
\end{equation}
where $\beta_3\equiv 0$ can be employed if $\lambda_0$ is not close to the zero-dispersion wavelength.

\section{Numerical methods for ultra-fast pulses in single-mode fibers}\label{sec:smssfm}

\subsection{Split-step Fourier approach}
In order to develop a numerical solution method, Eq.~(\ref{eq:pulse}) is commonly written in the form
\begin{equation}\label{eq:split} 
\frac{\partial A}{\partial z}= \underbrace{\left(-\frac{\alpha}{2}-
i\frac{\beta_2}{2}\frac{\partial^2}{\partial T^2}+\frac{\beta_3}{6}\frac{\partial^3}{\partial T^3}\right)}_{\cal D}A+
\underbrace{i\gamma\left(|A|^2 + i S \frac{1}{A}\frac{\partial}{\partial T}\left(A|A|^2\right)
-T_R \frac{\partial|A|^2}{\partial T} \right)}_{\cal N}A, 
\end{equation}
where we denote with ${\cal D}(A)$ the operator of all terms linear in $A$ and with ${\cal N}(A)$ the operator of all nonlinearities. Using these
definitions, we write Eq.~(\ref{eq:split}) in short as 
\begin{equation}\label{eq:splitbase}
\frac{\partial A}{\partial z} = ({\cal D}+{\cal N})A. 
\end{equation}
If one assumes ${\cal D}$ and ${\cal N}$ to be independent of $z$, Eq.~(\ref{eq:splitbase}) can be integrated exactly and the solution
at $z+h$ reads
\begin{equation}\label{eq:splitbase2}
A(z+h,T) = \exp(h({\cal D}+{\cal N}))A(z,T).
\end{equation}
The last expression forms the basis of split-step numerical methods \cite{Agrawal-07}. Note, however, that the operators ${\cal D}$ and 
$\cal N$ in general do not commute and that it corresponds to an $O(h)$ approximation to replace Eq.~(\ref{eq:splitbase2}) with
$\exp(h {\cal D})\exp(h {\cal N})A(z,T)$. A commonly used symmetric approximation is \cite{Strang-68,Glowinski-03}
\begin{equation}\label{eq:ssfm}
 A(z+h,T) = \exp\left(\frac{h}{2}{\cal D}\right)\exp(h{\cal N})\exp\left(\frac{h}{2}{\cal D}\right) A(z,T).
\end{equation}
Utilizing the Baker-Campbell-Hausdorff formula for expanding two non-commuting operators, Eq.~(\ref{eq:ssfm}) can be proven to be an $O(h^2)$ 
approximation \cite{Muslu-Erbay-03}. Comprehensive descriptions of the split-step approach for simulating pulse propagation in fibers are given for instance 
by Agrawal \cite{Agrawal-07} and Hohage \& Schmidt \cite{Hohage-02}. The efficiency of the SSFM, especially for longer propagation 
distances, as required for modeling optical communication lines, can be improved by taking solution adaptive steps in space as proposed 
by Sinkin {\em et al.} \cite{Sinkin-03}. 

Alternatively, one may also construct a fractional step splitting method by solving
\begin{equation} \frac{\partial A}{\partial z}= {\cal D} A\;,\quad \frac{\partial A}{\partial z}= {\cal N}(A)A = \bar {\cal N}(A) \end{equation}
successively, which we approximate with the symmetric fractional step method 
\begin{subequations}\label{eq:fullsplit}\begin{align}
A^\ast & = \exp\left(\frac{h}{2}{\cal D}\right) A(z,T)\;, \label{eq:fullsplita} \\
A^{\ast\ast} & = A^\ast+h\bar {\cal N}(A^\ast)\;, \label{eq:fullsplitb} \\
A(z+h,T) & = \exp\left(\frac{h}{2}{\cal D}\right) A^{\ast\ast}. \label{eq:fullsplitc}
\end{align}\end{subequations}
Note that step (\ref{eq:fullsplitb}) is written here as a simple explicit Euler method to motivate the fundamental idea but 
schemes described below are in fact more complicated. 

\subsection{Linear sub-steps}\label{sec:linear}
Since the dispersion parameters $\beta_2$ and $\beta_3$ are very small, discretization of the temporal derivatives in ${\cal D}$ by finite
differences and approximation in physical time is no viable option. Instead, Fourier transformation into frequency space is 
commonly applied. The linear operator then becomes
\begin{equation}\label{eq:fourier}
\exp\left(\frac{h}{2}{\cal D}\right)A(z,T) = {\cal F}^{-1}\exp\left[\frac{h}{2}\left(i\frac{\beta_2}{2}\omega^2-i\frac{\beta_3}{6}\omega^3-
\frac{\alpha}{2}\right)\right]{\cal F}A(z,T) ,
\end{equation}
where ${\cal F}$ and ${\cal F}^{-1}$ denote Fourier and inverse Fourier transformation, respectively. In the practical
implementation, discrete Fourier transformation needs to be used and for $\omega$ we employ the discrete frequency spectrum
\begin{equation} \left\{j\Delta\omega: j\in\mathbb{Z} \;\wedge\; -N\le j\le N-1\right\} \end{equation}
with spectral width $\Delta\omega=\pi/(N\Delta T)$. Here, it is assumed that the temporal window
traveling with the pulse is discretized with $2N$ points (note that discrete Fourier transformation algorithms 
are specially efficient if the number of points is a power of 2), $\Delta T$ denotes the temporal discretization width and
the temporal window has the extensions $\left[-N\Delta T,(N-1)\Delta T\right]$. 

\subsection{Nonlinear sub-steps}\label{sec:nonlinear}
The nonlinear operator $\cal N$ of the split-step method (\ref{eq:ssfm}) is discretized in physical space. Utilizing
$|A|^2=A\bar A$ to eliminate $1/A$, we write ${\cal N}(A)$ in the form
\begin{equation}\label{eq:ncentral} 
{\cal N}(A) = i\gamma\left(|A|^2 + i S \bar A\frac{\partial A}{\partial T}+ 
\left[iS-T_R\right] \frac{\partial|A|^2}{\partial T} \right).
\end{equation}
A consistent numerical method can be constructed by simply approximating the temporal derivatives in Eq.~(\ref{eq:ncentral})
by complex-valued first-order central differences and applying Eq.~(\ref{eq:ssfm}). The resulting split-step scheme would be second-order accurate in time and space. 
However, is is also clear that central finite differences will result in Gibbs phenomena (cf. \cite{Lax-06}) when strong self-steepening occurs 
or the propagation of an initially discontinuous signal needs to be simulated. 

An alternative approach for handling ${\cal N}(A)$ is to apply forward and inverse Fourier transformation individually to the 
derivatives, cf. \cite{Long-08}. For instance, in (\ref{eq:ncentral}) one simply replaces $\bar A\partial_T A$ and 
$\frac{\partial|A|^2}{\partial T}$ with $\bar A {\cal F}^{-1} (i\omega {\cal F} (A))$ and ${\cal F}^{-1} (i\omega {\cal F} (|A|^2)$, 
respectively, thereby neglecting the dependence of $\bar A$ on $T$. The result is class of numerical operators
that would generally not be consistent in the strict mathematical sense with ${\cal N}(A)$ and that are not uniquely defined, with
different authors arriving at slightly different disretizations of Eq.~(\ref{eq:ncentral}), cf. \cite{Long-08} and \cite{Amorim-09}.
Therefore, we have opted to pursue a different approach, which can handle self-steepening
and arbitrary signal shapes without artificial numerical oscillations. This method is based on solving
\begin{equation}\label{eq:pftrans} 
\frac{\partial A}{\partial z}= \underbrace{\left(-\frac{\alpha}{2}-
i\frac{\beta_2}{2}\frac{\partial^2}{\partial T^2}+\frac{\beta_3}{6}\frac{\partial^3}{\partial T^3}\right)}_{\cal D}A+
\underbrace{i\gamma\left(A|A|^2 + i S \frac{\partial}{\partial T}\left(A|A|^2\right)
-T_R A \frac{\partial|A|^2}{\partial T} \right)}_{\bar {\cal N}(A)}
\end{equation} 
within the fractional step method (\ref{eq:fullsplit}).
Specific to our approach is that we discretize and numerically solve the complete sub-operator
\begin{equation}\label{eq:nonlin} 
\frac{\partial A}{\partial z}= \bar {\cal N}(A) = i\gamma\left(A|A|^2 + 
i S \frac{\partial}{\partial T}\left(A|A|^2\right)-T_R A \frac{\partial|A|^2}{\partial T} \right)
\end{equation}
directly. Using the Madelung transformation \cite{Madelung-27,Spiegel-80} 
$A(z,t)=\sqrt{I(z,t)}e^{i\phi(z,t)}$, one can transform Eq.~(\ref{eq:nonlin}) into the
equivalent system of partial differential equations 
\begin{subequations} \begin{align} \frac{\partial I}{\partial z}+3 \gamma S I\frac{\partial I}{\partial T} = 0, \\
\frac{\partial \phi}{\partial z}+ \gamma S I\frac{\partial \phi}{\partial T} + 
\gamma T_R \frac{\partial I}{\partial T} = \gamma I \end{align} \end{subequations}
of the real-valued quantities intensity $I$ and phase $\phi$. If we write the latter in the form
\begin{equation}
\label{eq:nonlinadvec}
 \frac{\partial}{\partial z}\left[\begin{array}{c} I \\ \phi \end{array}\right] +
\left[\begin{array}{cc} 3 \gamma S I & 0 \\ \gamma T_R & \gamma S I \end{array}\right]
\frac{\partial}{\partial T}\left[\begin{array}{c} I \\ \phi \end{array}\right] = \left[\begin{array}{c} 0 \\ \gamma I \end{array}\right],
\end{equation}
its structure as a hyperbolic advection problem 
\begin{equation}\label{eq:advec}
 \frac{\partial\bf q}{\partial z}+{\bf M}({\bf q})\frac{\partial\bf q}{\partial T}={\bf s}({\bf q}) 
\end{equation}
with ${\bf q}=(I,\phi)^T$ becomes apparent. The matrix ${\bf M}({\bf q})$ has the eigenvalues $\lambda_1=3 \gamma S I$,
$\lambda_2=\gamma S I$ and a unique eigendecomposition for $I\ne 0$. Here we propose a
numerical method for (\ref{eq:nonlinadvec}) that considers the characteristic information, i.e., the sign of the eigenvalues 
of ${\bf M}({\bf q})$ for constructing one-sided (aka ``upwinded'') differences for the temporal derivatives, 
as it is required for a reliable and robust method following the theory of hyperbolic problems (cf. \cite{Smoller-82}).

Again, we adopt an operator splitting technique and, instead of discretizing (\ref{eq:advec}) directly, alternate between solving 
the homogeneous partial differential equation
\begin{equation}\label{eq:split1}
\partial_z{\bf q}+{\bf M}({\bf q})\partial_T{\bf q}=0
\end{equation}
and the ordinary differential equation
\begin{equation}\label{eq:split2}
\partial_z{\bf q} = {\bf s}({\bf q}) 
\end{equation}
successively, using the updated data from the preceding step as initial condition. A first-order accurate upwind scheme for (\ref{eq:split1}) 
can be derived easily based on the discrete update formula \cite{LeVeque-02}
\begin{equation}\label{eq:fd}
{\bf q}_j^{n+1} = {\bf q}_j^n - \frac{h}{\Delta T}\left({\bf\hat M}^-({\bf q}_{j+1},{\bf q}_j) \Delta{\bf q}_{j+1/2}^n + 
{\bf\hat M}^+({\bf q}_{j},{\bf q}_{j-1})\Delta{\bf q}_{j-1/2}^n\right) 
\end{equation}
with $\Delta{\bf q}_{j+1/2}^n = {\bf q}_{j+1}^n-{\bf q}_j^n$, where we assume a computational grid with equidistant mesh widths $\Delta T$ in time
indexed with $j\in\mathbb{Z}$, where $-N\le j\le N-1$, cf. Sect.~\ref{sec:linear}. The spatial update steps are 
indexed by $n\in\mathbb{N}_0$. In general, the matrices ${\bf \hat M}^+$ and ${\bf \hat M}^-$ indicate 
decompositions of ${\bf M}$ with only positive and negative eigenvalues, respectively. However, in the case of Eq.~(\ref{eq:nonlinadvec}) 
the eigenvalues have the same sign, which depends solely on the sign of $\gamma$ (since $I\ge 0$). Based on (\ref{eq:fd}), we
construct a straightforward upwind scheme for Eq.~(\ref{eq:nonlinadvec}) that reads 
\begin{subequations}\label{eq:upwindop1} \begin{align}\label{eq:upwindop1st1}
I_j^{n+1} & = I_j^n - \frac{h}{\Delta T}\left[3 \gamma S\tilde I^n_j\Delta I^n_j\right], \\ \label{eq:upwindop1st1b}
\bar\phi_j^{n+1} & = \phi_j^n - \frac{h}{\Delta T}\left[\gamma T_R \Delta I^n_j + \gamma S \tilde I^n_j \Delta \phi^n_j\right], \\
\phi_j^{n+1} & = \bar\phi_j^{n+1}+h\gamma I_j^{n+1} \label{eq:upwindop1st2}
\end{align}\end{subequations}
with 
\begin{eqnarray*}
\tilde I^n_j = \frac{1}{2}\left(I^n_j+I^n_{j-1}\right)\;,\quad \Delta I^n_j = I^n_j-I^n_{j-1} \qquad \text{for $\gamma>0$,} \\
\tilde I^n_j = \frac{1}{2}\left(I^n_j+I^n_{j+1}\right)\;,\quad \Delta I^n_j = I^n_{j+1}-I^n_j \qquad \text{for $\gamma<0$.} 
\end{eqnarray*}
When computing the phase difference $\Delta \phi^n_j$, it of crucial importance to remember that phase is given only modulo $2\pi$. 
Here, we have obtained reliable and stable results by ensuring that the smallest possible difference $\Delta \phi^n_j$ 
modulo $2\pi$ is applied in (\ref{eq:upwindop1st1b}). Using the auxiliary variable 
\begin{equation} \label{eq:upwindop1st5} 
\Delta\theta = \left\{ \begin{array}{ll}
\phi^n_j-\phi^n_{j-1}, & \text{for $\gamma>0$}, \\
\phi^n_{j+1}-\phi^n_j, & \text{for $\gamma<0$}.\\
\end{array} \right.\end{equation}
and $\Delta\tau_j^n=\min\left\{|\Delta\theta_j^n|,|\Delta\theta_j^n+2\pi|,|\Delta\theta_j^n-2\pi|\right\}$ we evaluate $\Delta \phi^n_j$ as
\begin{equation}\label{eq:phimod} 
\Delta\phi^n_j = \left\{ \begin{array}{ll}
\Delta\theta_j^n \,, & \text{if $|\Delta\theta_j^n|=\Delta\tau_j^n$}, \\
\Delta\theta_j^n+2\pi \,, & \text{if $|\Delta\theta_j^n+2\pi|=\Delta\tau_j^n$}, \\
\Delta\theta_j^n-2\pi \,, & \text{if $|\Delta\theta_j^n-2\pi|=\Delta\tau_j^n$}. \\
\end{array} \right.\end{equation}
The scheme (\ref{eq:upwindop1}) is of first-order accuracy and thereby entirely free
of producing numerical oscillations in the approximation of Eq.~(\ref{eq:pftrans}) provided that the stability condition
\begin{equation}\label{eq:cfl}
3|\gamma| S \underset{j}\max\left\{I^n_j\right\}\frac{h}{\Delta T}\le 1
\end{equation}
is satisfied. Our present implementation guarantees (\ref{eq:cfl}) under all circumstances by having the ability to adaptively take 
$k$ steps with step size $\Delta z$ with $h=k\Delta z$ within the central, nonlinear sub-step  
(\ref{eq:fullsplitb}) when required. Note, however, that for all computations presented in here
the stability conditions (\ref{eq:cfl}) was always already satisfied for $k=1$. 

To complete the algorithmic description we remark that we set $I_j^0:=|A_j^\ast|^2$ and $\phi^0_j:=\mathrm{arg}(A_j^\ast)$ after 
sub-step (\ref{eq:fullsplita}) and compute $A_j^{\ast\ast}=\sqrt{I_j^k} e^{i\phi^k_j}$ before step (\ref{eq:fullsplitc}). Periodic
boundary conditions could be implemented by one layer of halo points. But note that thanks to the
directional dependence, inherent to (\ref{eq:upwindop1}) and (\ref{eq:upwindop1st5}), it suffices to update only the 
upstream halo point, that is the one with index $j=-N-1$ for $\gamma>0$ and the one with $j=N$ in case $\gamma<0$ before applying
the upwind scheme.

\subsection{High-resolution upwind scheme}\label{sec:2ndorder}
To enable overall second-order numerical accuracy of the fractional step method (\ref{eq:fullsplit}), 
in case the solution is smooth and differentiable, it is necessary to extend the homogeneous nonlinear update (\ref{eq:fd}) to 
a {\em high-resolution scheme}. For this purpose, we have developed a special MUSCL-type slope-limiting technique of the solution
vector ${\bf q}$. Originally proposed by van Leer for hyperbolic equations in conservation law form \cite{vanLeer-79}, application to 
quasilinear systems is not apparent. Inspired by Ketcheson \& LeVeque \cite{Ketcheson-LeVeque-08}, we formulate our high-resolution 
method as  
\begin{equation}\label{eq:wpslope2}
{\bf q}_j^{n+1} \!\!=\! {\bf q}_j^n - \frac{h}{\Delta T}\!\left({\bf\hat M}^-\!\Delta{\bf q}_{j+1/2}^\star \!+\! 
{\bf\hat M}^+\!\Delta{\bf q}_{j-1/2}^\star\!+\!{\bf\hat M}\Delta{\bf q}_j^\star \!\right) 
\end{equation}
with $\Delta{\bf q}_{j+1/2}^\star={\bf q}_{j+1}^l-{\bf q}_j^r$, $\Delta{\bf q}_{j-1/2}^\star={\bf q}_j^l-{\bf q}_{j-1}^r$, and 
$\Delta {\bf q}_j^\star={\bf q}_j^r-{\bf q}_j^l$. Here, ${\bf q}_j^{l/r}$ refers to slope-limited values constructed for each
component of ${\bf q}$ separately as
\begin{equation}\label{eq:int} 
q^r_j = \bar q_j + \textstyle\frac{1}{4} \sigma_j, \qquad
q^l_j = \bar q_j - \textstyle\frac{1}{4} \sigma_j 
\end{equation}
with reconstructed linear local slope 
\begin{equation}\label{eq:slope}
\sigma_j = \Phi\left(\frac{\Delta_{j-\frac{1}{2}}}{\Delta_{j+\frac{1}{2}}}\right)\Delta_{j+\frac{1}{2}} + 
\Phi\left(\frac{\Delta_{j+\frac{1}{2}}}{\Delta_{j-\frac{1}{2}}}\right)\Delta_{j-\frac{1}{2}}
\end{equation}
with $\displaystyle \Delta_{j-1/2}=\bar q_j-\bar q_{j-1}$, $\displaystyle \Delta_{j+1/2}=\bar q_{j+1}-\bar q_j$.
In the latter, $\Phi(\cdot)$ is a typical limiter function, where we utilize in here exclusively 
the {\em van Albada} limiter 
\begin{equation}\label{eq:vanalbada}
\Phi(r)=\max\left(0,(r^2+r)/(1+r^2)\right).
\end{equation}
To permit second-order accuracy overall, we do not utilize in (\ref{eq:int}) 
the discrete values from the previous step ${\bf q}^n$ but instead intermediate values ${\bf\bar q}$ computed as 
\begin{equation}\label{eq:wpslope1}
{\bf\bar q}_j = {\bf q}_j^n - \frac{h}{2 \Delta T}\left({\bf\hat M}^-\Delta{\bf q}_{j+1/2}^n + {\bf\hat M}^+\Delta{\bf q}_{j-1/2}^n\right) .
\end{equation}
The consecutive application of (\ref{eq:wpslope2}) and (\ref{eq:wpslope1}) corresponds to an explicit 2-step Runge-Kutta method
in the spatial update. Finally, a second-order accurate symmetric operator splitting \cite{Strang-68,Glowinski-03} is employed to 
integrate Eq.~(\ref{eq:split2}) before and after the high-resolution scheme. Thanks to the simplicity of ${\bf s}({\bf q})$ 
using an explicit Euler method for this step is equivalent to an explicit 2-step Runge-Kutta update.

We want to point out that the first-order method (\ref{eq:upwindop1}) as well as
the MUSCL-based second-order scheme are equally applicable for $T_R=0$ and especially in the singular case $S=0$, 
which allows deactivation of Raman scattering and/or self-steepening if desired. Note that for $S=0$ or 
$\underset{j}\max\left\{I^n_j\right\}=0$, the stability condition (\ref{eq:cfl}) is trivially satisfied.

\subsection{Simulation of a propagating pulse}\label{sec:pulse1}
In order to demonstrate the described numerical method we simulate the propagation of a Gaussian pulse with initial shape
\begin{equation}\label{eq:gaussian} A(0,T) = \sqrt{P_0}\exp\left(-\frac{1+iC}{2}\frac{T^2}{T_0^2}\right) \end{equation}
in a homogeneous fiber. The fiber is assumed to be lossless ($\alpha=0$) for simplicity as the omitted linear weakening of the signal is 
unproblematic for any numerical scheme. We use the SSFM in line with Eq.~(\ref{eq:fullsplit}) with second-order accurate 
upwind-based nonlinear operator, cf. Sect.~\ref{sec:2ndorder}, and Van Albada slope-limiter (\ref{eq:vanalbada}). 

Used parameters correspond to a typical ultra-short communication pulse with $P_0=0.625\,\mathrm{mW}$, 
$T_0=80\,\mathrm{fs}$, and no chirp, i.e. $C=0$. The central wavelength is set to $\lambda_0=1550\,\mathrm{nm}$, from which 
one computes the self-steepening parameter 
$S=\lambda_0/2\pi c$, with $c$ denoting the speed of light in vacuum. Raman scattering is activated with $T_R=3\,\mathrm{fs}$. 
Realistic fiber parameters $\beta_2=0.5\,\mathrm{ps^2 km^{-1}}$, $\beta_3=0.07\,\mathrm{ps^3 km^{-1}}$ and $\gamma=0.1\,\mathrm{W/m}$ are used.  
For this configuration, the second-order dispersion length is just $L_d=T_0^2/|\beta_2|\approx 12.8\,\mathrm{m}$, the third-order dispersion length is 
$T_0^3/|\beta_3|\approx 7.31\,\mathrm{m}$, and the nonlinear length is $L_{nl}=\left(\gamma P_0\right)^{-1}=16\,\mathrm{km}$. 
The pulse is assumed to travel a distance of just $L_{\max}=1\,\mathrm{km}$ and the simulated temporal window moving with the pulse has the width 
$[-30\,\mathrm{ps},30\,\mathrm{ps}-\Delta T]$. 

\begin{figure}[t]
\sidecaption[t]
\includegraphics[angle=-90,origin=c,width=7cm]{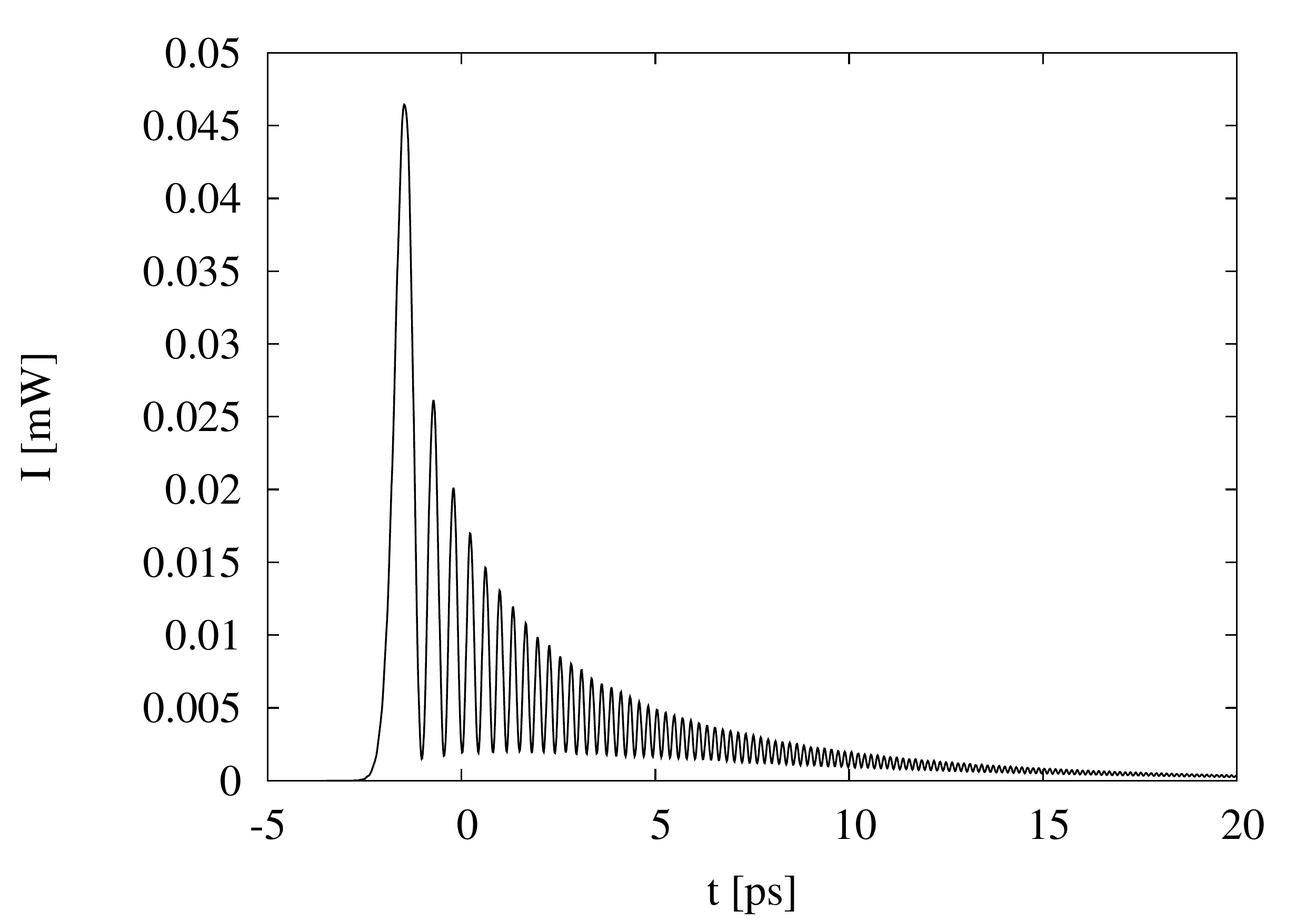}\vspace*{-1cm}
\caption{Simulated signal at $L_{\max}=1\,\mathrm{km}$ (temporal window enlarged) for Benchmark 1. The initially Gaussian pulse, 
cf. Eq.~(\ref{eq:gaussian}), with half width 
$T_0=80\,\mathrm{fs}$ has broadened severely because of second-order dispersion. Asymmetric high-frequency oscillations 
have been added by third-order dispersion effects. Maximal signal strength is reduced by a factor of $\sim 13.9$.}
\label{pic:bench1}
\end{figure}

Figure~\ref{pic:bench1} shows the computed solution using a temporal discretization of 
$2N$ points for $N=2048$ and after taking $M=100$ spatial steps of equal size of $h=10\,\mathrm{m}$. Because of the very small
second- and third-order dispersion lengths, typical for ultra-fast pulses, the final signal shape is clearly dominated by
dispersion effects. Second-order dispersion has introduced severe pulse broadening, reducing the maximum in power by a factor of $\sim 13.9$; 
third-order dispersion has added high-frequency oscillations. 

\begin{figure}[t]
\parbox[t]{5.5cm}{
\centering\includegraphics[angle=-90,width=5.75cm]{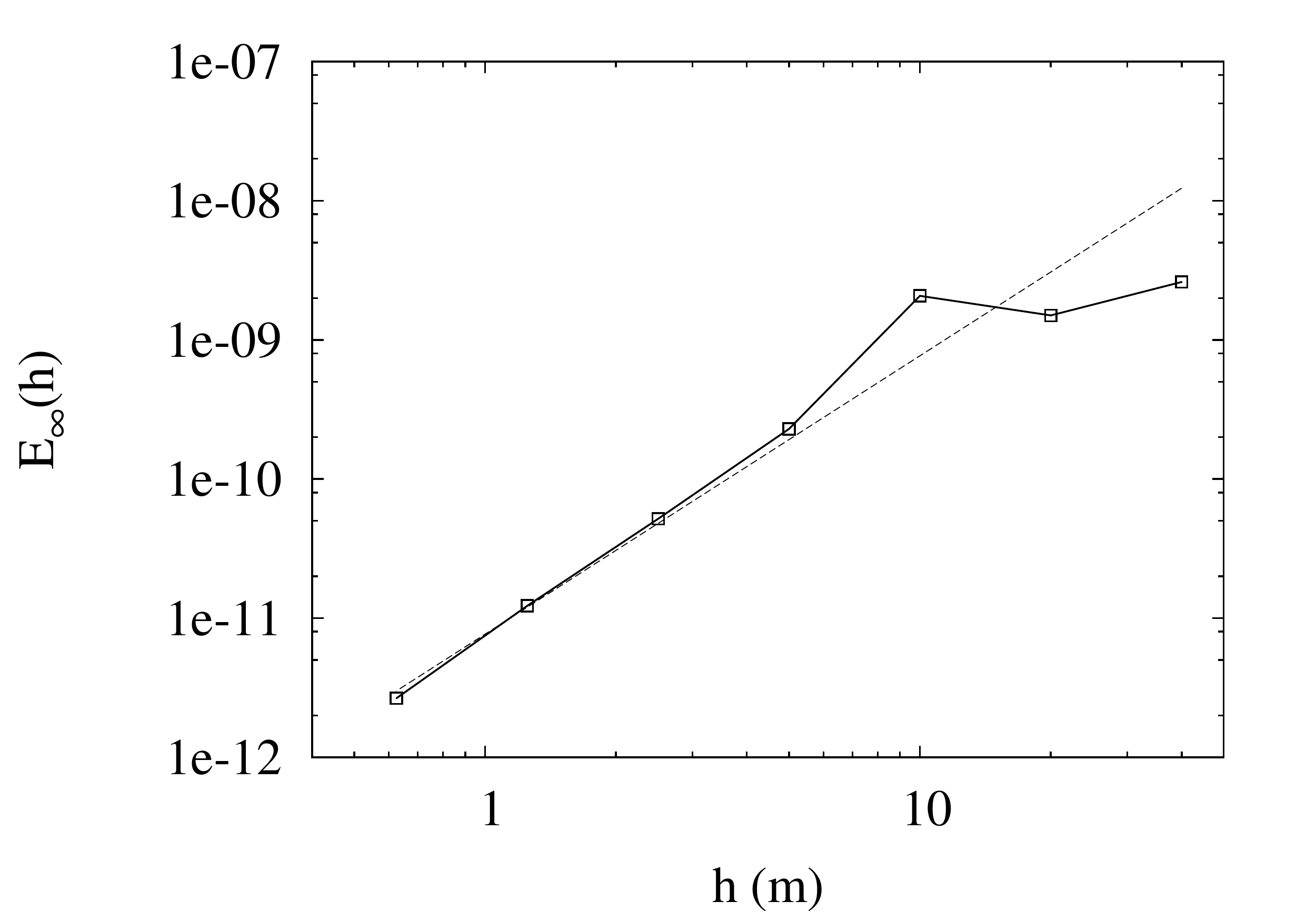}

\caption{Numerical error $E_\infty$ over $h$ for Benchmark 1. The dotted line corresponds to ideal second order approximation accuracy.} 
\label{pic:bench1b}}\hspace*{0.5cm}\parbox[t]{5.5cm}{\centering\includegraphics[angle=-90,width=5.75cm]{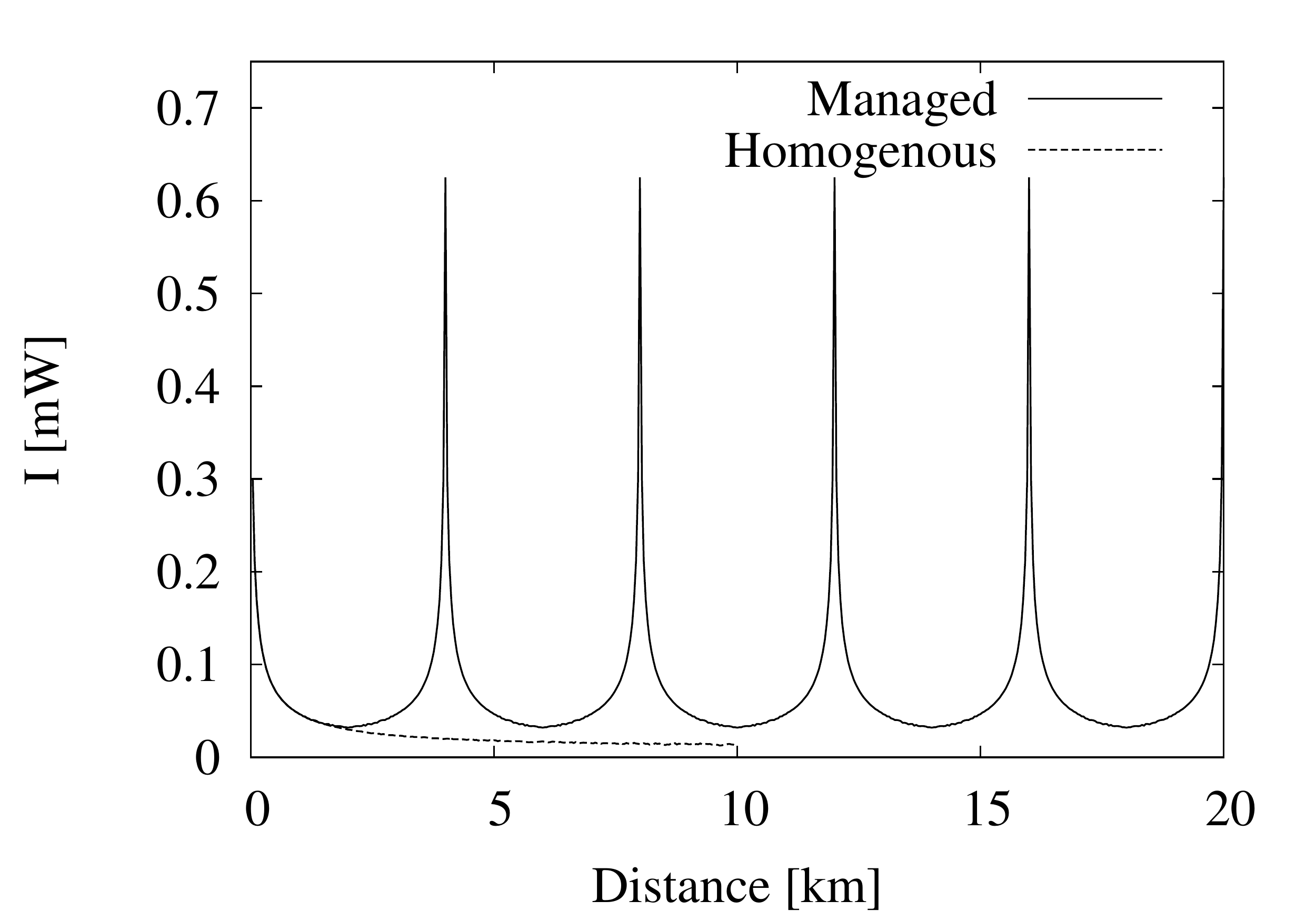}

\caption{Benchmark 2: Maximal power over distance with and without dispersion management.} 
\label{pic:bench2a}}
\end{figure}

A detailed numerical analysis verifies the convergence and expected order of accuracy of the scheme. Starting from $N=512$ and 
$h=40\,\mathrm{m}$ ($M=25$ steps), in each successive computation the number of Fourier modes and spatial steps is doubled. 
The numerical error at $L_{\max}$ is measured for the intensity of the signal in the discrete maximum norm 
\begin{equation}\label{eq:error}
E_\infty=\underset{j\in\{-N,N-1\}}\max|I_j-I^\mathrm{ref}(j\Delta T)|,
\end{equation}
where a highly resolved result with $N=131,072$ and $M=6400$ is used as reference solution $I^\mathrm{ref}$. 
Figure~\ref{pic:bench1b} visualizes the numerical error $E_\infty$ over $h$ and it is eminent that the method achieves almost 
perfect second-order approximation accuracy and reliable, robust convergence. A more detailed numerical study of the second-order accurate 
upwind-based SSFM including comparisons with several alternative numerical methods can be found in \cite{Deiterding-etal-13}.

\subsection{Spatially dependent fiber parameters}\label{sec:dispmanage}
Continued propagation of the pulse of Fig.~\ref{pic:bench1} will invariably lead to a signal which has broadened to such an extent
that it can not be used for digital communication. Yet, this problem can be compensated surprisingly easily by combining fiber sections
with positive and negative dispersion characteristics into a single communication line. This technique is called {\em dispersion
management} and has been studied extensively both theoretically and numerically because of its practical significance for long-distance
fiber optical communication \cite{Malomed-97,Richardson-00,Atre-07}. Instead of Eq.~(\ref{eq:pftrans}), one considers the extended variant 
\begin{multline}\label{eq:splitz} 
\frac{\partial A}{\partial z}= \left(-\frac{\alpha(z)}{2}-
i\frac{\beta_2(z)}{2}\frac{\partial^2}{\partial T^2}+\frac{\beta_3(z)}{6}\frac{\partial^3}{\partial T^3}\right) A \\
+i\gamma(z)\left(A|A|^2 + i S\frac{\partial}{\partial T}\left(A|A|^2\right)
-T_R A \frac{\partial|A|^2}{\partial T} \right)
\end{multline}
as governing equation. Adopting the practical viewpoint that the spatial numerical steps of any SSFM will be significantly larger than the 
spatial extension corresponding to the used temporal simulation window moving with the pulse, a straightforward numerical method for 
Eq.~(\ref{eq:splitz}) can be constructed by simply averaging the spatially dependent parameters between discrete propagation steps, i.e. by using 
\begin{equation} \bar\beta_{\left\{2,3\right\},j}=\frac{2}{h}\int\limits_{z_j}^{z_j+\frac{h}{2}} \beta_{\left\{2,3\right\}}(\xi)d\xi\;,\quad 
\bar\alpha_j=\frac{2}{h}\int\limits_{z_j}^{z_j+\frac{h}{2}} \alpha(\xi)d\xi\end{equation}
in the linear numerical operator (\ref{eq:fourier}) and by using
\begin{equation} \bar\gamma_j=\frac{1}{h}\int\limits_{z_j}^{z_j+h} \gamma(\xi)d\xi \end{equation}
in the nonlinear operator approximating (\ref{eq:nonlin}). 

\begin{figure}[t]
\centering\includegraphics[angle=-90,width=6cm]{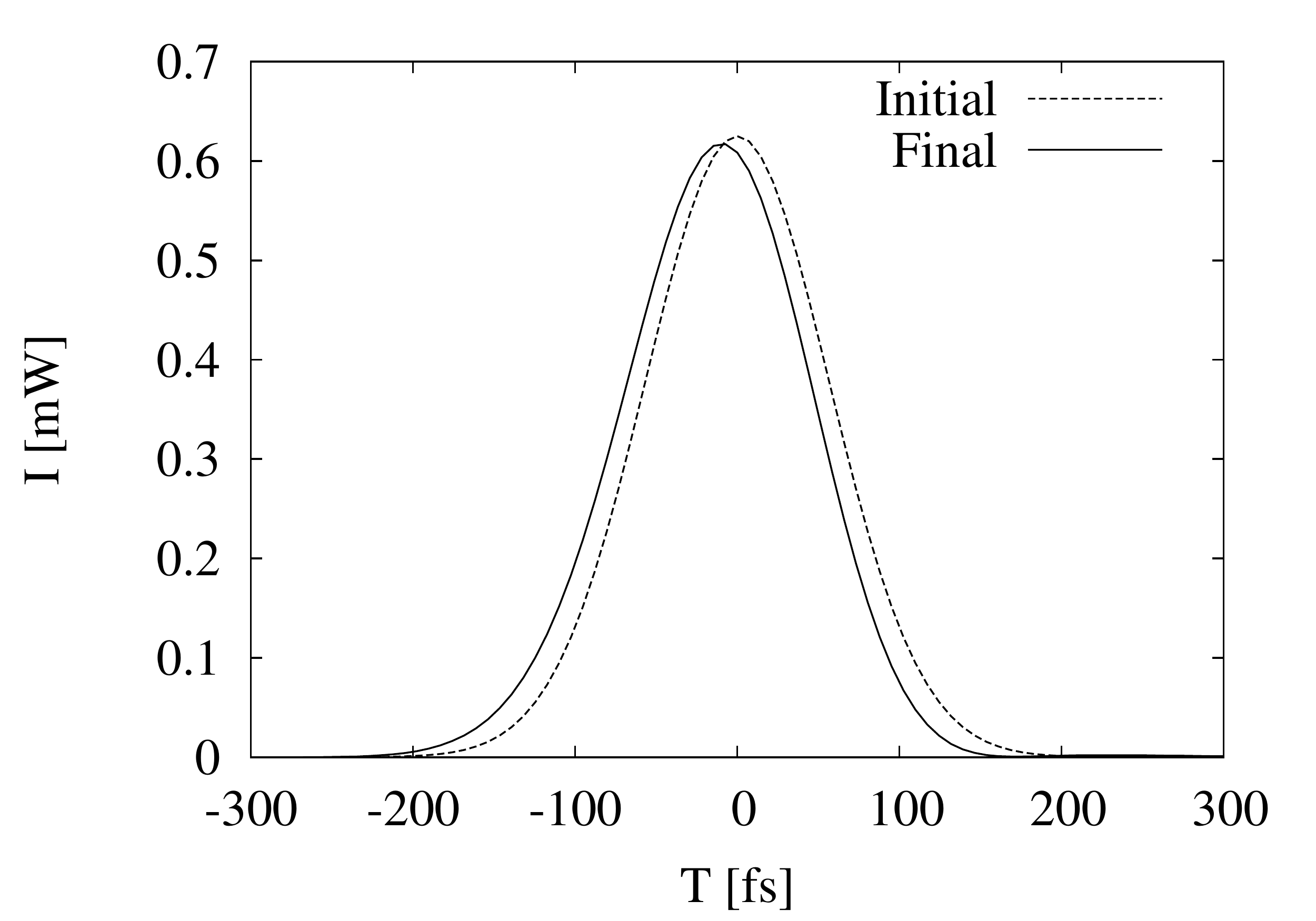}\includegraphics[angle=-90,width=6cm]{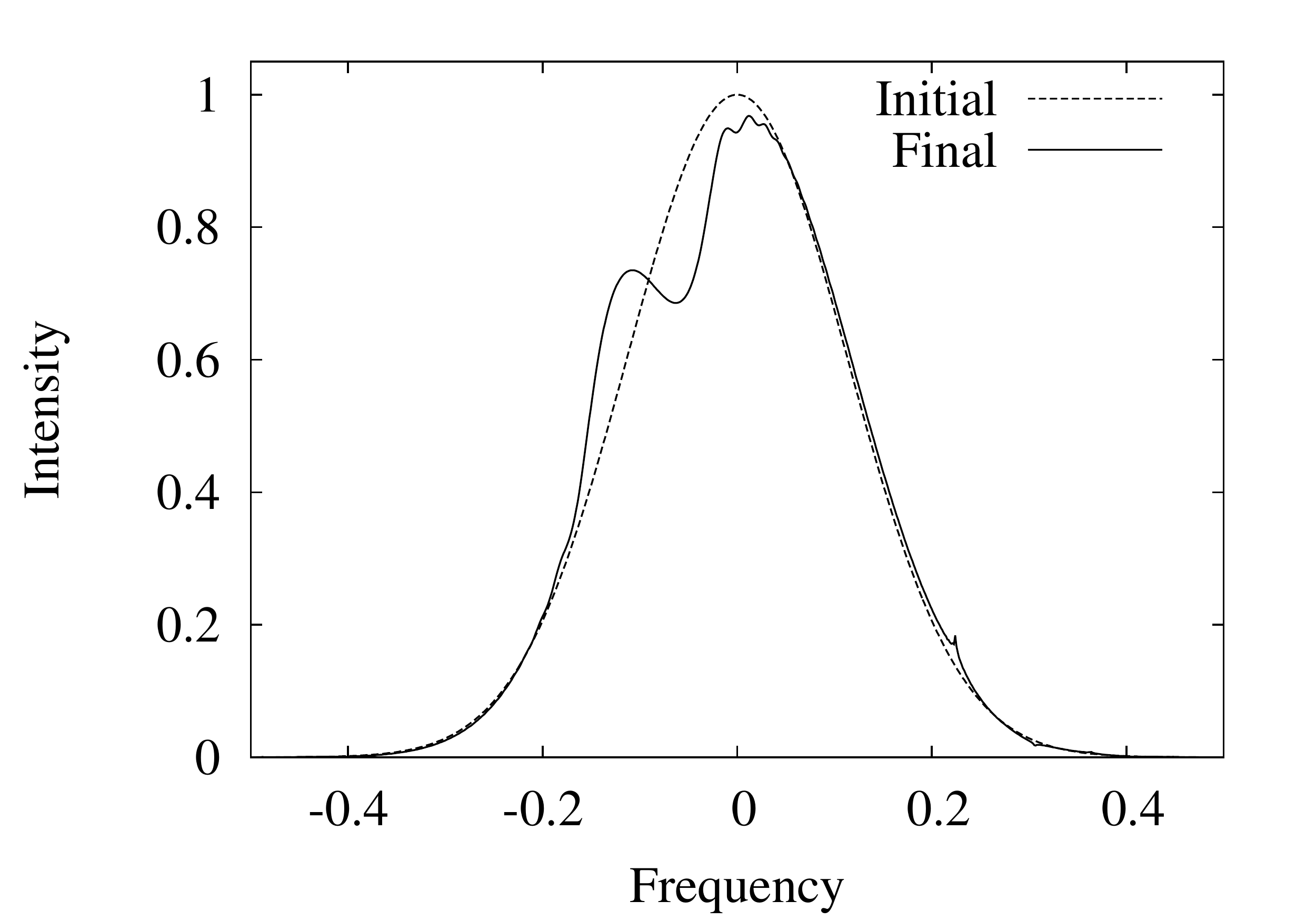}

\caption{Benchmark 2: Pulse shape and spectrum after propagating $100\,\mathrm{km}$ or experiencing 25 soliton-like oscillations from alternating signs of 
dispersion parameters.} 
\label{pic:bench2b} 
\end{figure}

In practice, very sophisticated dispersion management designs might be 
employed (for instance, Guo \& Huang \cite{Guo-04} propose an exponential decrease of $|\beta_2|$ to accommodate better for linear loss). 
Here, we simply extend the example of Sect.~\ref{sec:pulse1} and alternate the sign of $\beta_2$ and $\beta_3$ every $2\,\mathrm{km}$. 
All other parameters are unaltered and for an example computation we use $N=4096$ and $h=40\,\mathrm{m}$ ($M=2500$) to simulate a pulse propagation
over a distance of $100\,\mathrm{km}$. In the fiber sections with negative dispersion parameters, the pulse deterioration is effectively
reversed and the pulse shape mostly recovered. The pulse is undergoing a soliton-like oscillation with a period of $4\,\mathrm{km}$, which can
be inferred from Fig.~\ref{pic:bench2a}. This graphic compares the pulse power peak over distance in the simulation with periodic
dispersion management and when the computation of the previous section is continued to a length of $10\,\mathrm{km}$. In Fig.~\ref{pic:bench2b}
are compared the shape and spectra of the initial Gaussian pulse and of the signal after propagating for $100\,\mathrm{km}$. The observed
slight signal delay and spectral modification is the combined effects of the nonlinearities. If $\gamma=0$ is used, the initial signal
is exactly recovered.

\section{Governing equations for two interacting ultra-fast pulses}\label{sec:twomode}
Data throughput can be increased significantly if multiple optical fields of different wavelengths propagate simultaneously inside
the fiber. However, these fields would interact with one another through all the fiber nonlinearities. Additionally if three or more 
fields are initially present, even new signal fields can be induced (aka {\em four-wave mixing} \cite{Agrawal-07}). Therefore, we consider in the following
only the case of two interacting signal fields propagating through an optical fiber, for which there is already some agreement about
the structure of the governing equations in the literature \cite{Kalithasan-08}. Extensions of the ultra-fast pulse propagation equation (\ref{eq:pulsefast}) 
to three or more interacting fields are still a topic of active research.

We assume two pulses at carrier frequencies $\omega_0^{(1)}$, $\omega_0^{(2)}$, and two nonlinear constants $\gamma_1$, $\gamma_2$. It is further assumed that 
the cross-phase modulation of each frequency can be expressed for all higher order nonlinear terms by positive factors $B_1$, $B_2$, the cross-phase
modulation in the Kerr nonlinearity by factors $C_1$, $C_2$. Extending Eq.~(\ref{eq:pulsefast}) accordingly, we use the model equations
\begin{subequations}\label{eq:tmodeul} \begin{align}
\displaystyle\frac{\partial A_1}{\partial z} =&\displaystyle -\frac{\alpha_1}{2}A_1-\beta_1^{(1)}\frac{\partial A_1}{\partial t}-
i\frac{\beta_2^{(1)}}{2}\frac{\partial^2A_1}{\partial t^2}+\frac{\beta_3^{(1)}}{6}\frac{\partial^3A_1}{\partial t^3}+
\displaystyle i\gamma_1 \left(|A_1|^2+C_1|A_2|^2\right)A_1 \nonumber \\
& \displaystyle -\frac{\gamma_1}{\omega_0^{(1)}}\left[\frac{\partial\left(|A_1|^2A_1\right)}{\partial t}+
B_1\frac{\partial\left(|A_2|^2A_1\right)}{\partial t}\right]-
i\gamma_1T_R\left[\frac{\partial|A_1|^2}{\partial t}+B_1\frac{\partial|A_2|^2}{\partial t}\right]A_1, \\ 
\displaystyle\frac{\partial A_2}{\partial z} =&\displaystyle -\frac{\alpha_2}{2}A_2-\beta_1^{(2)}\frac{\partial A_2}{\partial t}-
i\frac{\beta_2^{(2)}}{2}\frac{\partial^2A_2}{\partial t^2}+\frac{\beta_3^{(2)}}{6}\frac{\partial^3A_2}{\partial t^3}+
\displaystyle i\gamma_2\left(|A_2|^2+C_2|A_1|^2\right)A_2 \nonumber \\
& \displaystyle -\frac{\gamma_2}{\omega_0^{(2)}}\left[\frac{\partial\left(|A_2|^2A_2\right)}{\partial t}+
B_2\frac{\partial\left(|A_1|^2A_2\right)}{\partial t}\right]-
i\gamma_2T_R\left[\frac{\partial|A_2|^2}{\partial t}+B_2\frac{\partial|A_1|^2}{\partial t}\right]A_2 .  
\end{align}\end{subequations}
Note that (\ref{eq:tmodeul}) encompasses the model actually adopted for simulation by Kalithasan {\em et al.} in \cite{Kalithasan-08}. 
Using $\beta_1^{(j)}=1/v_j$ and the transformation $T\equiv t-\beta_1^{(j)} z$ into retarded time yields
\begin{subequations}\label{eq:tmpframe} \begin{align}
\displaystyle\frac{\partial A_1}{\partial z} =&\displaystyle -\frac{\alpha_1}{2}A_1\qquad\quad
-i\frac{\beta_2^{(1)}}{2}\frac{\partial^2A_1}{\partial T^2}+\frac{\beta_3^{(1)}}{6}\frac{\partial^3A_1}{\partial T^3}+
\displaystyle i\gamma_1 \left(|A_1|^2+C_1|A_2|^2\right)A_1 \nonumber \\
& \displaystyle -\gamma_1 S_1 \left[\frac{\partial\left(|A_1|^2A_1\right)}{\partial T}+
B_1\frac{\partial\left(|A_2|^2A_1\right)}{\partial T}\right]-
i\gamma_1T_R\left[\frac{\partial|A_1|^2}{\partial T}+B_1\frac{\partial|A_2|^2}{\partial T}\right]A_1,  \\ 
\displaystyle\frac{\partial A_2}{\partial z} =&\displaystyle -\frac{\alpha_2}{2}A_2-\delta\frac{\partial A_2}{\partial T}-
i\frac{\beta_2^{(2)}}{2}\frac{\partial^2A_2}{\partial T^2}+\frac{\beta_3^{(2)}}{6}\frac{\partial^3A_2}{\partial T^3}+
\displaystyle i\gamma_2\left(|A_2|^2+C_2|A_1|^2\right)A_2 \nonumber \\
& \displaystyle -\gamma_2 S_2\left[\frac{\partial\left(|A_2|^2A_2\right)}{\partial T}+
B_2\frac{\partial\left(|A_1|^2A_2\right)}{\partial T}\right]-
i\gamma_2T_R\left[\frac{\partial|A_2|^2}{\partial T}+B_2\frac{\partial|A_1|^2}{\partial T}\right]A_2 ,  
\end{align}\end{subequations}
with $\delta=(v_1-v_2)/(v_1v_2)$ representing the group velocity mismatch between both fields.
As before we use $S_k=1/\omega_0^{(k)}$ for $k=1,2$.

In the regime of pico-second pulses, that is for pulses with $T_0\gg 1\,\mathrm{ps}$, two-mode extensions of Eq.~(\ref{eq:pulseslow}) are well 
established. Setting $S_k= 0$, $T_R=0$ and using $C_k=2$ in (\ref{eq:tmpframe}), 
we obtain the frequently used \cite{Agrawal-07} cross-phase modulation model 
\begin{subequations}\label{eq:tmd} \begin{align}
\displaystyle\frac{\partial A_1}{\partial z} =& \displaystyle\underbrace{\left(-\frac{\alpha_1}{2}\qquad\quad-
i\frac{\beta_2^{(1)}}{2}\frac{\partial^2}{\partial T^2}+\frac{\beta_3^{(1)}}{6}\frac{\partial^3}{\partial T^3}\right)}_{{\cal D}^{(1)}}A_1 +
& \displaystyle \underbrace{i\gamma_1 \left(|A_1|^2+2|A_2|^2\right)}_{{\cal N}^{(1)}}A_1,\\
\displaystyle\frac{\partial A_2}{\partial z} =& \displaystyle\underbrace{\left(-\frac{\alpha_2}{2}- \delta\frac{\partial}{\partial T}-
i\frac{\beta_2^{(2)}}{2}\frac{\partial^2}{\partial T^2}+\frac{\beta_3^{(2)}}{6}\frac{\partial^3}{\partial T^3}\right)}_{{\cal D}^{(2)}}A_2 +
& \displaystyle \underbrace{i\gamma_2 \left(|A_2|^2+2|A_1|^2\right)}_{{\cal N}^{(2)}}A_2,
\end{align}\end{subequations}
which we write as 
\begin{equation} \frac{\partial A_1}{\partial z}=\left({\cal D}^{(1)}+{\cal N}^{(1)}(A_1,A_2)\right)A_1,\quad 
\frac{\partial A_2}{\partial z}=\left({\cal D}^{(2)}+{\cal N}^{(2)}(A_1,A_2)\right)A_2. 
\end{equation}

\section{Numerical methods for two interacting ultra-fast pulses}\label{sec:smsstm}

\subsection{Extended split-step Fourier method}\label{sec:exssfm}
Taking advantage of the fact that the linear operators ${\cal D}^{(k)}$ only need to be applied to each field $A_k$, 
a SSFM for approximating solutions of system (\ref{eq:tmd}) -- in line with Eq.~(\ref{eq:ssfm}) -- is easily constructed as 
\begin{subequations}\begin{align}\label{eq:pstmssfm1}
A_1^\ast & = \exp\left(\frac{h}{2}{\cal D}^{(1)}\right)A_1, & A_2^\ast & = \exp\left(\frac{h}{2}{\cal D}^{(2)}\right)A_2, \\
\label{eq:pstmssfm2}
A_1^{\ast\ast} & = \exp\left(h {\cal N}^{(1)}(A_1^\ast,A_2^\ast)\right)A_1^\ast,
&  A_2^{\ast\ast} & = \exp\left(h {\cal N}^{(2)}(A_1^{\ast\ast},A_2^\ast)\right)A_2^\ast \\
\label{eq:pstmssfm3}
A_1(z+h) & = \exp\left(\frac{h}{2}{\cal D}^{(1)}\right)A_1^{\ast\ast},
& \!\!\!\!\!\!\!\!\!\!\!\!\!\!A_2(z+h) & = \exp\left(\frac{h}{2}{\cal D}^{(2)}\right)A_2^{\ast\ast}.
\end{align}\end{subequations}
Obviously, the numerical operators of (\ref{eq:pstmssfm2}) and (\ref{eq:pstmssfm3}) acting on each fields can be executed consecutively.
The linear operator ${\cal D}^{(1)}$ is identical to (\ref{eq:fourier}). For ${\cal D}^{(2)}$ we have 
\begin{equation}\label{eq:fourier2m} 
\exp\left(\frac{h}{2}{\cal D}^{(2)}\right)A_2 = {\cal F}^{-1}\exp\left[\frac{h}{2}\left(-i\delta\omega+i\frac{\beta^{(2)}_2}{2}\omega^2-i\frac{\beta^{(2)}_3}{6}\omega^3-
\frac{\alpha_2}{2}\right)\right]{\cal F}A_2.
\end{equation}
A second-order accurate scheme can be expected if (\ref{eq:pstmssfm2}) is replaced with a symmetric splitting scheme such as
\begin{subequations}\label{eq:pstmssfm4}\begin{align} 
A_1^\star & = \exp\left(\frac{h}{2} {\cal N}^{(1)}(A_1^\ast,A_2^\ast)\right)A_1^\ast,\\
A_2^{\ast\ast} & = \exp\left(h {\cal N}^{(2)}(A_1^\star,A_2^\ast)\right)A_2^\ast, \\
A_1^{\ast\ast} & = \exp\left(\frac{h}{2} {\cal N}^{(1)}(A_1^\star,A_2^{\ast\ast})\right)A_1^\star.
\end{align}\end{subequations}

\subsection{Nonlinear sub-steps}\label{sec:tmnonlinear}
While the derivation of a SSFM for the simplified system (\ref{eq:tmd}) is apparently a straightforward task, formulation
of a reliable numerical method for the system of propagation equations for two coupled ultra-fast pulsed signals, (\ref{eq:tmpframe}), 
is more involved. In particular, when the equations of (\ref{eq:tmpframe}) are written in the form $\partial_z A_k = ({\cal D}^{(k)}+{\cal N}^{(k)})A_k$ one quickly finds
that due to the cross-phase coupling the factor $1/A_k$ of the self-steepening term cannot be eliminated from ${\cal N}^{(k)}$ as it was done to 
obtain Eq.~(\ref{eq:ncentral}). This leaves a singularity in the operator for vanishing signals and neither the centered difference
method nor particularly an ad hoc Fourier transformation technique, sketched both in the beginning of Sect.~\ref{sec:nonlinear},
are available anymore for numerical method construction. However, we will demonstrate subsequently how our upwind-based discretization technique 
of Sect.~\ref{sec:nonlinear} can be easily extended to (\ref{eq:tmpframe}), yielding a reliable and robust numerical method.

We start the derivation of the method by inserting the linear operators from (\ref{eq:tmd}) into (\ref{eq:tmpframe}) to obtain
\begin{equation}\label{eq:Dtmp} \begin{array}{lll}
\displaystyle\frac{\partial A_k}{\partial z} & = \displaystyle {\cal D}^{(k)}A_k+i\gamma_k\left(|A_k|^2+C_k|A_l|^2\right)A_k\\[1.5ex] 
& \displaystyle -\gamma_k S_k\left[\frac{\partial\left(|A_k|^2A_k\right)}{\partial T}+
B_k\frac{\partial\left(|A_l|^2A_k\right)}{\partial T}\right]
-i\gamma_kT_R\left[\frac{\partial|A_k|^2}{\partial T}+B_k\frac{\partial|A_l|^2}{\partial T}\right]A_k
\end{array}\end{equation}
for $k,l\in\{1,2\}$ and $k\ne l$.
In analogy to Sect.~\ref{sec:nonlinear}, we assume a fractional step approach in the spirit of Eq.~(\ref{eq:fullsplit})
that considers the linear operators with the update steps (\ref{eq:pstmssfm1}) and (\ref{eq:pstmssfm3}) and approximates the
nonlinear sub-operator equations
\begin{equation}\label{eq:tmnlop} \begin{array}{lll}
\displaystyle\frac{\partial A_k}{\partial z} =& \displaystyle\bar {\cal N}^{(k)}(A_k,A_l) = i\gamma_k\left(|A_k|^2+C_k|A_l|^2\right)A_k \\[1.5ex]
& \displaystyle -\gamma_k S_k\left[\frac{\partial\left(|A_k|^2A_k\right)}{\partial T}+
B_k\frac{\partial\left(|A_l|^2A_k\right)}{\partial T}\right]
-i\gamma_kT_R\left[\frac{\partial|A_k|^2}{\partial T}+B_k\frac{\partial|A_l|^2}{\partial T}\right]A_k. \end{array}\end{equation}
Using again Madelung transformation for each field, i.e. $A_k(z,t)=\sqrt{I_k(z,t)}e^{i\phi_k(z,t)}$, we obtain the 
transport equations for the intensities $I_k$ and the phases $\phi_k$ instead of (\ref{eq:tmnlop}) as
\begin{subequations}\label{eq:tmn}\begin{align}
 \frac{\partial I_k}{\partial z}+\gamma_k S_k\left[\left(3I_k+B_kI_l\right)\frac{\partial I_k}{\partial T}+
   2B_kI_k\frac{\partial I_l}{\partial T}\right] & =0, \\
\frac{\partial \phi_k}{\partial z}+\gamma_k S_k \left(I_k+B_kI_l\right)\frac{\partial \phi_k}{\partial T}
+\gamma_kT_R\left[\frac{\partial I_k}{\partial T}+B_k\frac{\partial I_l}{\partial T}\right] & =\gamma_k\left(I_k+C_kI_l\right).
\end{align}\end{subequations}
The latter defines a single system of advection equations that couples the fields $A_k$ and $A_l$.
Using the state vector $\mathbf{u}=\left(I_1,\phi_1,I_2,\phi_2\right)^T$, this system reads 
\begin{equation}\label{eq:tmhyp}
\frac{\partial\mathbf{u}}{\partial z}+ \mathbf{B}(\mathbf{u}) \frac{\partial\mathbf{u}}{\partial T} =\mathbf{r}(\mathbf{u}), 
\end{equation}
with matrix 
\begin{equation} \mathbf{B}(\mathbf{u})=\left[\begin{array}{cccc}
\gamma_1 S_1(3I_1+B_1I_2) & 0 & 2\gamma_1S_1B_1I_1 & 0 \\
\gamma_1 T_R & \gamma_1S_1(I_1+B_1I_2) & \gamma_1T_R B_1 & 0 \\
2\gamma_2 S_2 B_2 I_2 & 0 & \gamma_2S_2 (3I_2+B_2 I_1) & 0 \\
\gamma_2 T_R B_2 & 0 & \gamma_2 T_R & \gamma_2 S_2 (I_2+B_2 I_1) \end{array}\right] \end{equation}
and right hand side
\begin{equation} 
\mathbf{r}(\mathbf{u})= \left(0,\gamma_1\left(I_1+C_1I_2\right),0,\gamma_2\left(I_2+C_2I_1\right)\right)^T. 
\end{equation}
In order to verify the hyperbolicity of Eq.~(\ref{eq:tmhyp}) and for constructing an upwind scheme, one would require the eigendecomposition 
$\mathbf{B}=\mathbf{R} \Lambda \mathbf{R}^{-1}$. However, the necessary linear algebra is very involved and is additionally complicated
by the singular cases $I_k=0$, which have to be considered separately in order to construct a generally robust numerical scheme. To simplify
the latter, we have opted to use a splitting approach and update the fields $A_k$ and $A_l$ successively. Instead of solving the combined
system (\ref{eq:tmhyp}) we construct an approximation to (\ref{eq:tmn}) under the assumption that $I_l$ is independent
of $z$. Proceeding then as in Sect.~\ref{sec:nonlinear}, we write (\ref{eq:tmn}) as the advection system
\begin{equation}\label{eq:tmredhyp}
\frac{\partial\mathbf{q^{(k)}}}{\partial z}+ \mathbf{M^{(k)}}(\mathbf{q^{(k)}}) \frac{\partial\mathbf{q^{(k)}}}{\partial T} = \mathbf{s^{(k)}}(\mathbf{q^{(k)}}), 
\end{equation}
with vector of state $\mathbf{q}^{(k)}=(I_k,\phi_k,I_l)^T$, matrix
\begin{equation}  \mathbf{M^{(k)}}(\mathbf{q^{(k)}}) =\left[\begin{array}{ccc}
\gamma_k S_k ( 3 I_k+B_kI_l) & 0 & 2\gamma_kS_kB_kI_k \\ 
\gamma_k T_R & \gamma_k S_k (I_k+B_k I_l) & \gamma_k T_RB_k \\
0 & 0 & 0 \end{array}\right] \end{equation}
and source term
\begin{equation} \mathbf{s^{(k)}}(\mathbf{q^{(k)}}) = \left(0, \gamma_k\left(I_k+C_kI_l\right), 0 \right)^T. \end{equation}
The non-zero eigenvalues of $\mathbf{M^{(k)}}$ are $\gamma_k S_k ( 3 I_k+B_kI_l)$ and $\gamma_k S_k (I_k+B_k I_l)$. Since $B_k\ge 0$ and $I_{k/l}\ge 0$
hold true, both eigenvalues have again the same sign, solely determined by the sign of $\gamma_k$. Following the upwind approach again
we construct a first-order accurate method for (\ref{eq:tmredhyp}) as  
\begin{subequations}\label{eq:tm1st}\begin{align}\label{eq:tm1stI}
I_{k,j}^{n+1} & = I_{k,j}^n - \frac{h}{\Delta T}\gamma_k S_k \left[(3\tilde I^n_{k,j}+B_k\tilde I^n_{j,l}) \Delta I^n_{k,j}+
2 B_k \tilde I^n_{k,j}\Delta I^n_{l,j}\right], \\
\bar\phi_{k,j}^{n+1} & = \phi_{k,j}^n - \frac{h}{\Delta T}\gamma_k \left[T_R (\Delta I^n_{k,j}+B_k\Delta I^n_{l,j}) + 
S_k (\tilde I^n_{k,j}+B_k\tilde I^n_{l,j}) \Delta \phi^n_{k,j}\right], \\ \label{eq:tm1stphi}
\phi_{k,j}^{n+1} & = \bar\phi_{k,j}^{n+1}+h\gamma_k\left(I_{k,j}^{n+1} +C_k I_{l,j}^n\right), 
\end{align}\end{subequations}
where 
\begin{eqnarray*}
\tilde I^n_{k/l,j}= \frac{1}{2}\left(I^n_{k/l,j}+I^n_{k/l,j-1}\right),\quad \Delta I^n_{k/l,j}= I^n_{k/l,j}-I^n_{k/l,j-1} \qquad \text{for $\gamma_k>0$,} \\
\tilde I^n_{k/l,j}= \frac{1}{2}\left(I^n_{k/l,j}+I^n_{k/l,j+1}\right),\quad \Delta I^n_{k/l,j}= I^n_{k/l,j+1}-I^n_{k/l,j},\qquad \text{for $\gamma_k<0$.} 
\end{eqnarray*}
As before, $\Delta\phi^n_{k,j}$ is evaluated modulo $2\pi$ using Eqs.~(\ref{eq:upwindop1st5}) and (\ref{eq:phimod}) and the stability condition
reads 
\begin{equation}\label{eq:cfl2}
 |\gamma_k| S_k \underset{j}\max\left\{3 I_{k,j}+B_kI_{l,j}\right\}\frac{h}{\Delta T}\le 1.
\end{equation}
By construction, the single-field upwind method (\ref{eq:tm1st}) computes only new values for $I_k$ and $\phi_k$, while
the intensity of the other field, $I_l$, is assumed to remain unchanged. In order to achieve an update of both fields, and thereby approximation of 
(\ref{eq:tmhyp}), we apply the single-field upwind scheme within a symmetric fractional-step splitting method, i.e. 
\begin{subequations}\label{eq:pstmssfm5}\begin{align} A_1^\star & = A_1^\ast + \frac{h}{2}\bar {\cal N}^{(1)}(A_1^\ast,A_2^\ast) ,\\
A_2^{\ast\ast} & = A_2^\ast + h \bar {\cal N}^{(2)}(A_1^\star,A_2^\ast), \\
A_1^{\ast\ast} & = A_1^\star+\frac{h}{2} \bar {\cal N}^{(1)}(A_1^\star,A_2^{\ast\ast}).
\end{align}\end{subequations}
A symmetric SSFM is obtained by applying expressions (\ref{eq:pstmssfm1}), (\ref{eq:pstmssfm5}) and (\ref{eq:pstmssfm3}) after another. 
Finally, the high-resolution technique, described in Sect.~\ref{sec:2ndorder}, is adopted to implement a second-order
accurate approximation to $\bar {\cal N}^{(k)}$, where we presently apply slope-limited reconstruction to $I_k$ and $\phi_k$ but not to $I_l$.

%\begin{eqnarray}
%\bar I_{k,j}^{n+1} & = & I_{k,j}^n - \frac{h}{\Delta T}\gamma_k S_k \left[(3\tilde I^n_{k,j}+B_k\tilde I^n_{j,l}) \Delta I^n_{k,j}\right] \\
%\bar\phi_{k,j}^{n+1} & = & \phi_{k,j}^n - \frac{h}{\Delta T}\gamma_k \left[T_R \Delta I^n_{k,j} + 
%S_k (\tilde I^n_{k,j}+B_k\tilde I^n_{l,j}) \Delta \phi^n_{k,j}\right], \\
%I_{k,j}^{n+1} & = & \bar I_{k,j}^{n+1} - \frac{2h}{\Delta T}\gamma_k S_k B_k \bar I^{n+1}_{k,j}\Delta I^n_{k,l} \\
%\phi_{k,j}^{n+1} & = & \bar\phi_{k,j}^{n+1}+h\gamma_k\left(\bar I_{k,j}^{n+1} +C_k I_{l,j}^n\right) - \frac{h}{\Delta T}\gamma_k T_R B_k\Delta I^n_{k,l}
%\end{eqnarray}

\subsection{Simulation of two interacting propagating pulses}\label{sec:pulse2}
We use a configuration with very strong nonlinearity and thereby nonlinear pulse interaction to assess the reliability of the 
derived two-mode method. A fiber without linear loss and third-order dispersion is assumed, i.e. $\alpha_{1,2}=0$ and $\beta_3^{(1,2)}=0$,
and Raman scattering is also deactivated by setting $T_R=0$. To enforce a strong influence of the nonlinearities we use 
$\beta_2^{(1,2)}=4\times 10^{-5}\,\mathrm{ps^2 km^{-1}}$, $\gamma_1=1\,\mathrm{W/m}$, and $\gamma_2=1.2\,\mathrm{W/m}$. Two 
unchirped pulses in the range of ultra-short communication pulses with $T_0^{(1,2)}=80\,\mathrm{fs}$ and power levels of 
$P_0^{(1)}=0.625\,\mathrm{mW}$ and $P_0^{(2)}=0.3125\,\mathrm{mW}$ are used. The first central wavelength is set to 
$\lambda_0^{(1)}=1550\,\mathrm{nm}$ and the second to $\lambda_0^{(2)}=1300\,\mathrm{nm}$.  The group velocity mismatch
parameter is set to $\delta = 0.015625\,\mathrm{fs/m}$ and the cross-phase modulation parameters read $B_{1,2}=C_{1,2}=2$. 

\begin{figure}[t]
\centering\includegraphics[angle=-90,width=6cm]{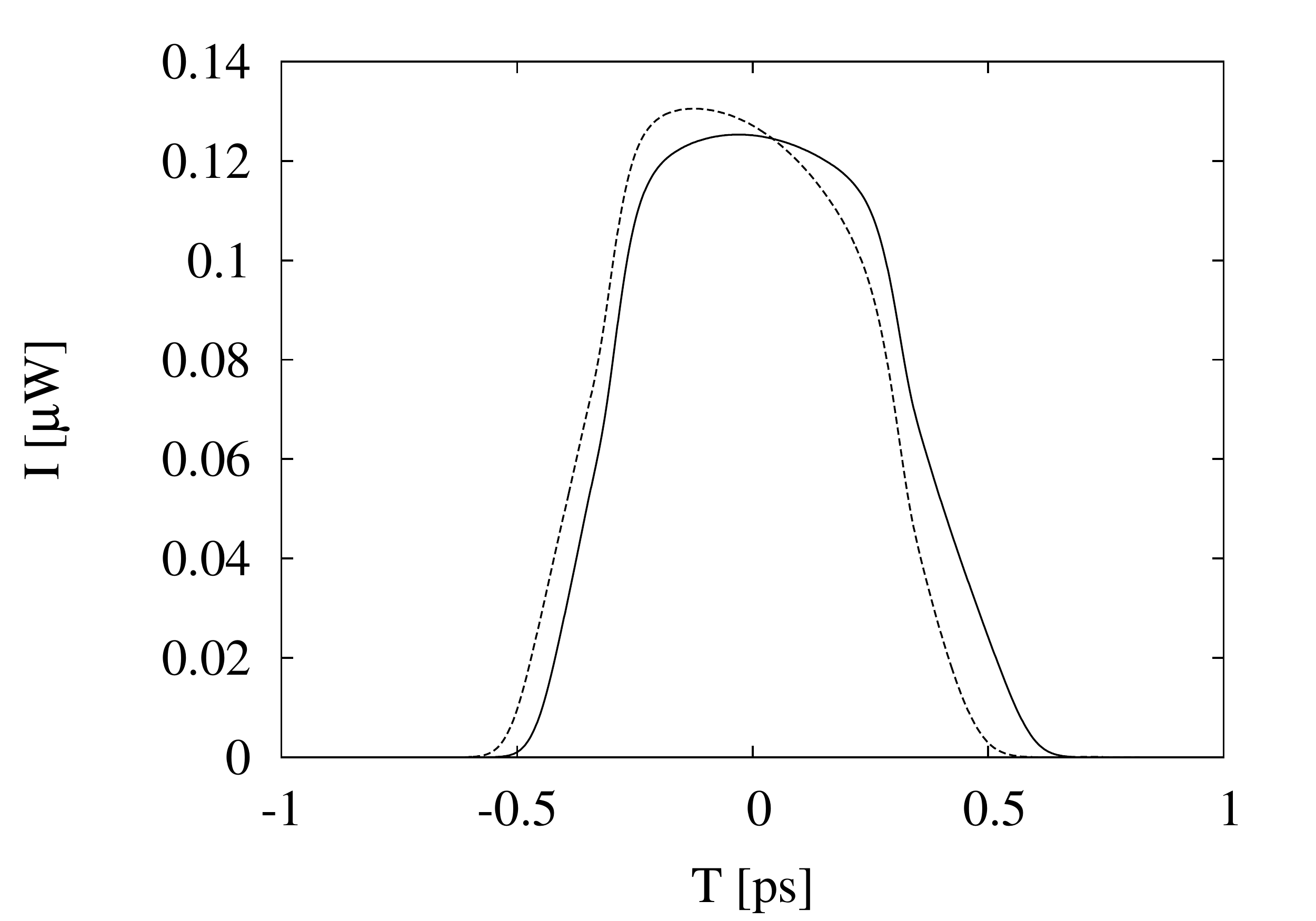}\includegraphics[angle=-90,width=6cm]{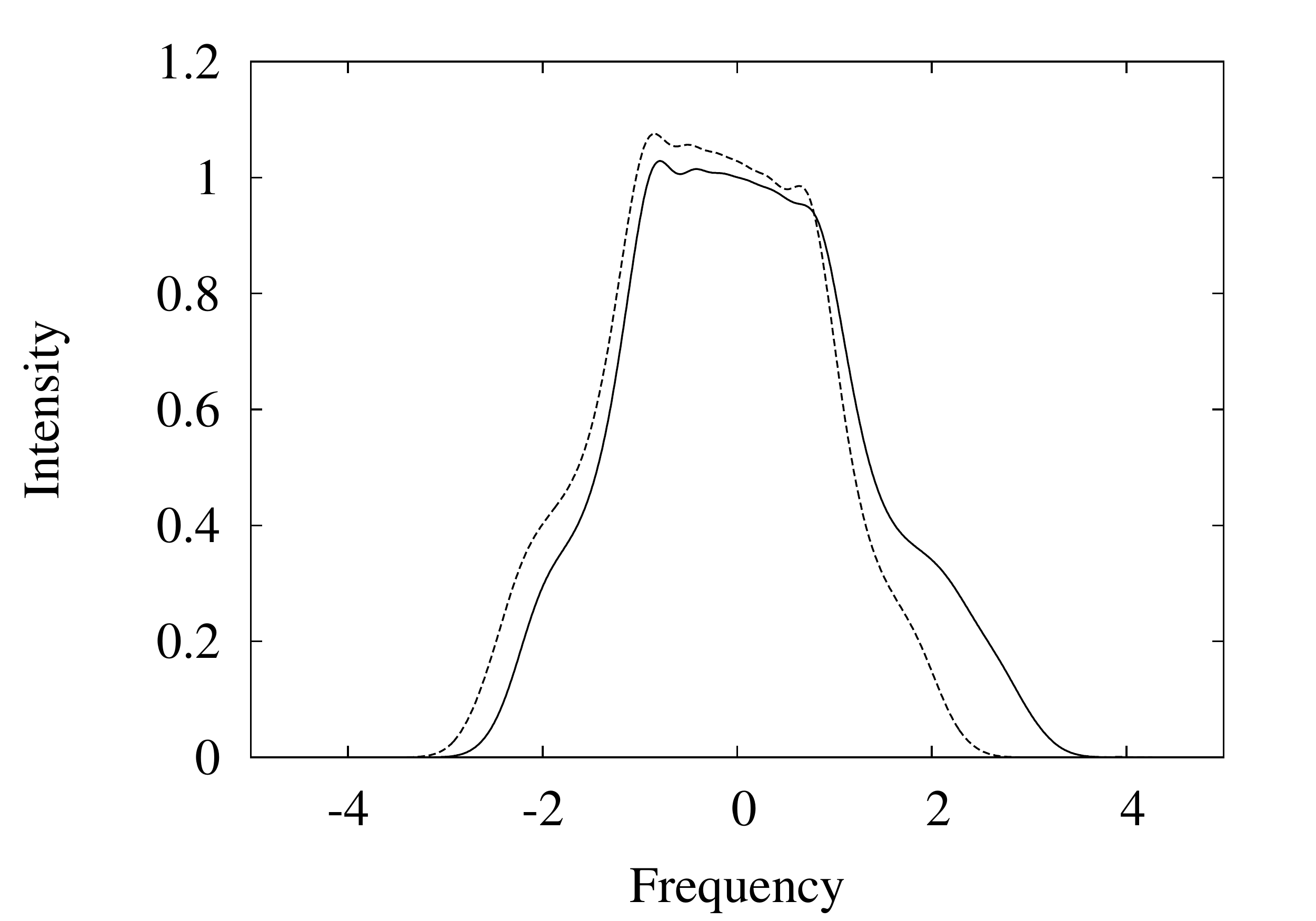}

\includegraphics[angle=-90,width=6cm]{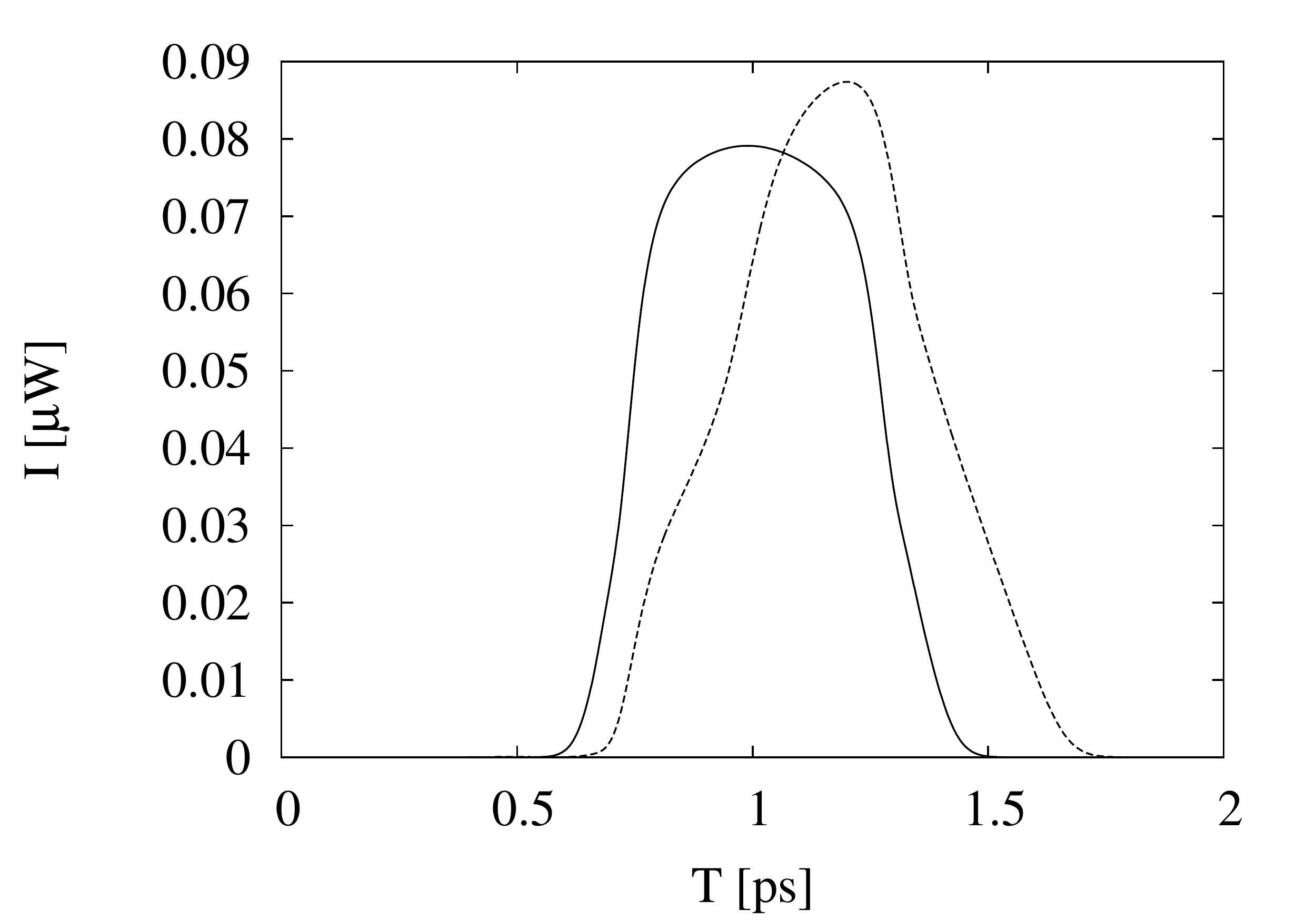}\includegraphics[angle=-90,width=6cm]{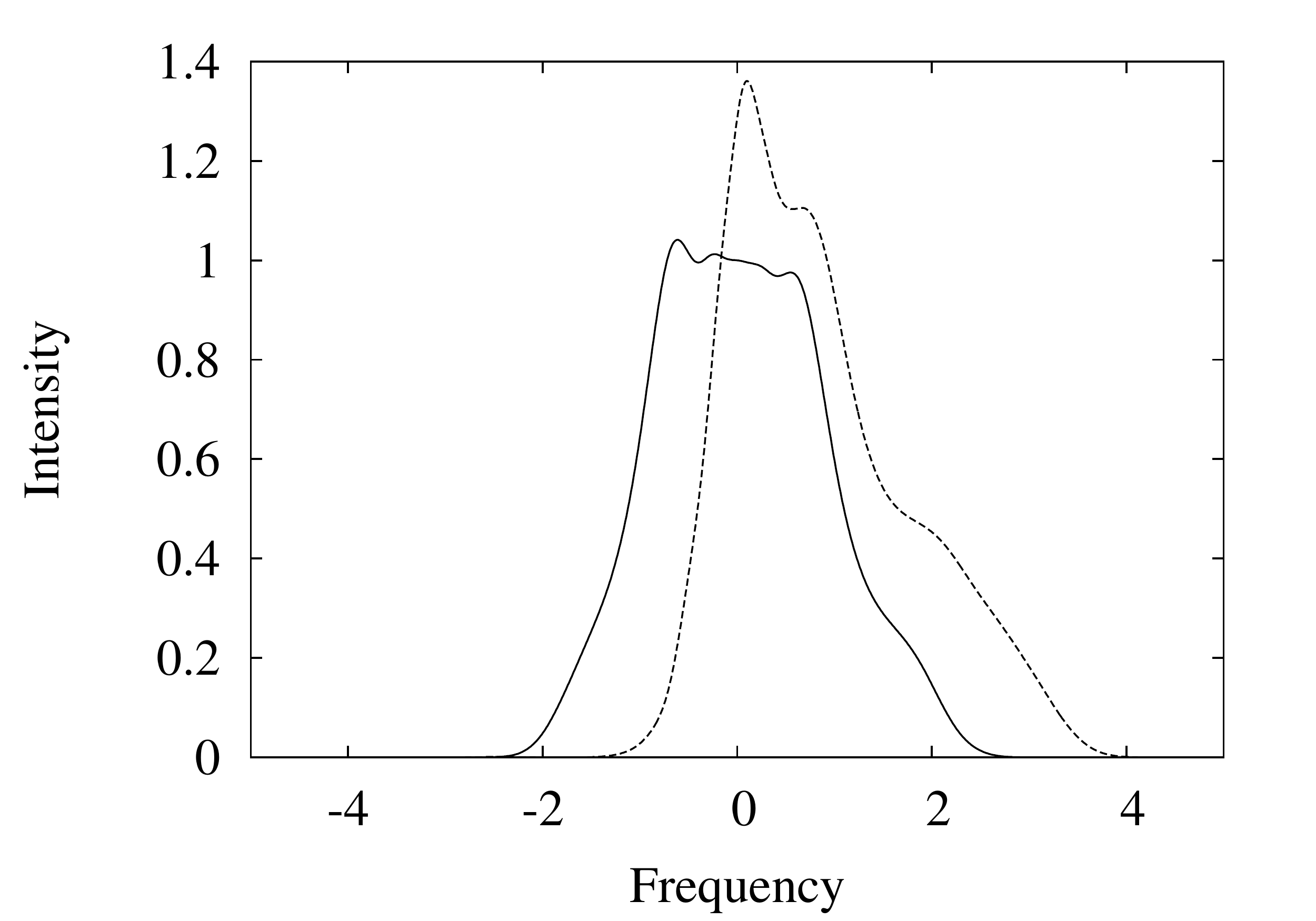}

\caption{Benchmark 3: Pulse shape and spectrum after propagating $64\,\mathrm{km}$ 
of two individual highly nonlinear ultra-short single-mode pulses (solid lines) 
and when the two pulses are interacting with one another in a two-mode fiber. The upper row corresponds to Pulse 1 with $P_0^{(1)}=0.625\,\mathrm{mW}$;
the lower row to Pulse 2 with $P_0^{(2)}=0.3125\,\mathrm{mW}$ and $\delta = 0.015625\,\mathrm{fs/m}$ causing the pulse to arrive $1\,\mathrm{ps}$ 
earlier.} 
\label{pic:bench3} 
\end{figure}

For this configuration, the second-order dispersion length is $L_d=160\,\mathrm{km}$ and the nonlinear lengths 
$L_{nl}^{(1)}=1.6\,\mathrm{km}$ and $L_{nl}^{(2)}=2.667\,\mathrm{km}$, respectively. The approximate optical shock distances 
\cite{Agrawal-07}, $z_s^{(1,2)}=\sqrt{e} L_{nl}^{(1,2)}\omega_0^{(1,2)} T_0/(3\sqrt{2})$ are $\sim 60.491\,\mathrm{km}$ and 
$\sim~\!\!120.216\,\mathrm{km}$, respectively. We use a propagation distance of $L_{\max}=64\,\mathrm{km}$, yielding a temporal shift
of the second pulse by exactly $1\,\mathrm{ps}$, and the temporal window has the width $[-4\,\mathrm{ps},4\,\mathrm{ps}-\Delta T]$.

In Fig.~\ref{pic:bench3} is shown the computed solution using a temporal discretization of $2N$ points for $N=2048$ and after taking $M=3200$ 
spatial steps of equal size of $h=20\,\mathrm{m}$. Additionally are shown the solutions if each pulse travels individually. These solutions
are computed by keeping all other parameters unchanged while setting $P_0^{(2)}\equiv 0$ and $P_0^{(1)} \equiv 0$, respectively. If only a single
field is present, our two-mode SSFM is identical to the previously developed single-mode SSFM, which was confirmed to be second-order 
accurate in Sect.~\ref{sec:pulse1}. Note that the single-mode solution of Pulse 1 was also used as a detailed computational benchmark in 
\cite{Deiterding-etal-13} and is thereby available as a reference. From Fig.~\ref{pic:bench3} it can be seen that the two 
non-interacting single-mode pulses exhibit a very similar shape and spectrum. However, in the two-mode model particularly the faster and 
weaker second pulse is significantly altered. Pulse 2, visualized in the lower row of Fig.~\ref{pic:bench3}, experiences considerable signal 
steepening from cross-phase modulation, which can be inferred especially from its spectrum.  

\begin{figure}[t]
\sidecaption[t]
\includegraphics[angle=-90,origin=c,width=7.2cm]{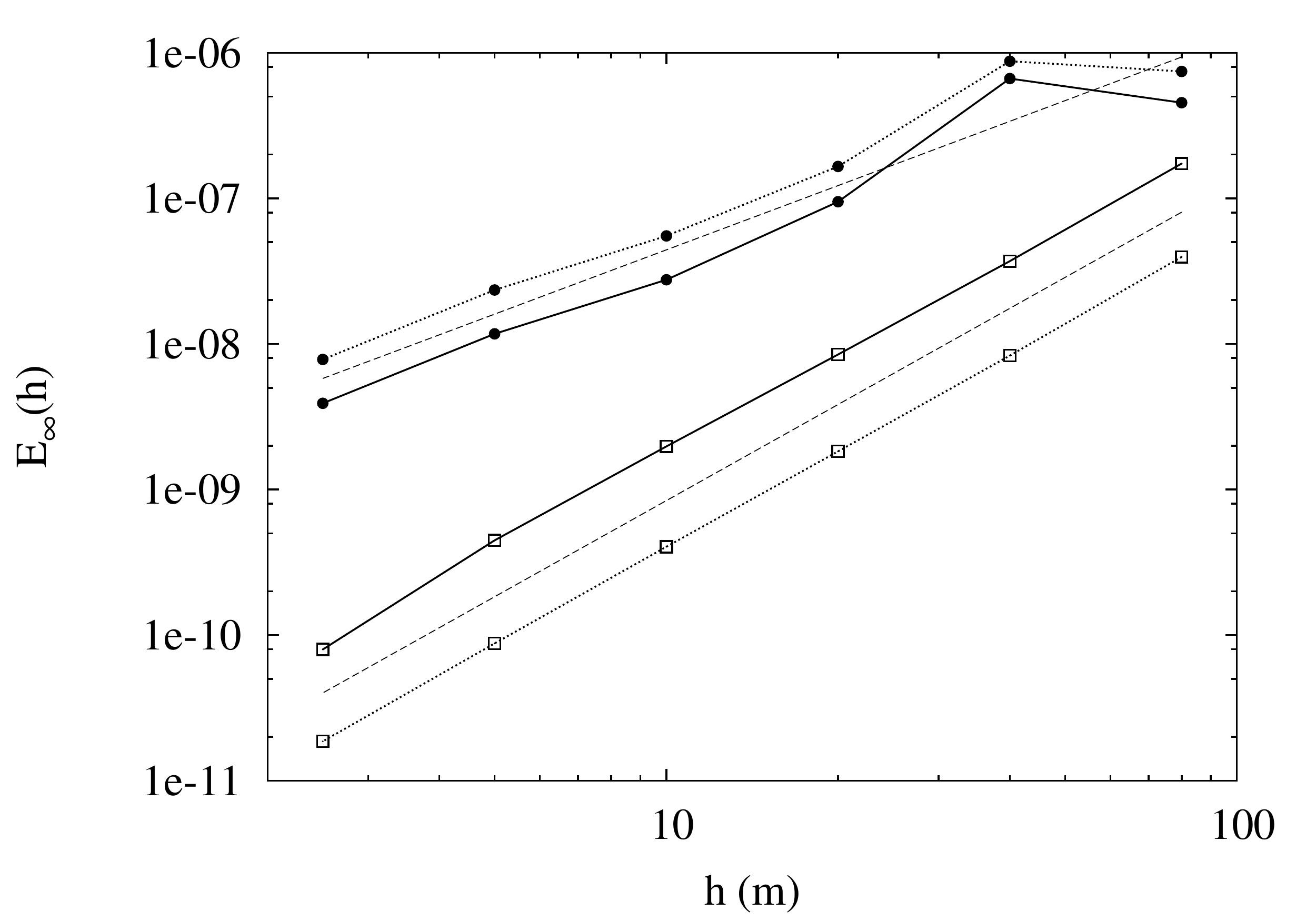}\vspace*{-1cm}
\caption{Numerical error $E_\infty$ over $h$ for Benchmark 3.  The respective error of Pulse 1 is marked with solid 
lines, the respective error of Pulse 2 with dotted lines. Single-pulse simulation results are indicated with open 
squares, the fully coupled simulations are marked with closed circles. The upper broken line corresponds 
to an order of accuracy $\sim 1.47$, the lower one to an order of accuracy of $\sim 2.20$.} 
\label{pic:bench3b} 
\end{figure}

We use the same technical approach as in Sect.~\ref{sec:pulse1} to quantify the numerical error and order of accuracy of the two-mode SSFM.
We double the temporal resolution consecutively starting from $N=512$ up to $N=16,384$ and simultaneously divide the spatial step size by a factor 
of 2 respectively, starting with $h=80\,\mathrm{m}$ ($M=800$ steps). The numerical error at $L_{\max}$ is measured for $I_{(1,2)}$ in the maximum norm, 
cf. Eq.~(\ref{eq:error}), where results computed with $N=32,768$ and $h=1.25\,\mathrm{m}$ ($M=51,200$) are used as respective reference solutions.
The computational errors of a series of fully coupled two-mode results as well as the errors of single-mode computations (cf. 
Fig.~\ref{pic:bench3}) of both individual pulses are plotted in Fig.~\ref{pic:bench3b}. In general, the example confirms that the proposed 
two-mode SSFM converges reliably and robustly even for a highly nonlinear coupled problem and performs identical beside round-off errors 
to the single-mode method of Sect.~\ref{sec:2ndorder} for uncoupled individual pulses. While the single-mode SSFM 
with limiter (\ref{eq:vanalbada}) actually achieves slight super-convergence in this test case (the measured order of accuracy is $\sim 2.20$),
the two-mode SSFM of Sect.~\ref{sec:tmnonlinear} with same limiter yields an approximate order of accuracy of $\sim 1.47$. 
One might attribute this behavior to the fractional step splitting treatment of the nonlinear operator, (\ref{eq:pstmssfm5}), 
however increasing the number of spatial steps up to a factor of 8 to possible reduce the splitting error of the nonlinear sub-operator
resulted only in marginally smaller numerical errors for this test case. 

\subsection{Spatially dependent fiber parameters}
As final test case, the coupled propagation of the two Gaussian pulses of the previous benchmark through the dispersion-managed communication 
line of Sect.~\ref{sec:dispmanage} is considered. Like in Sect.~\ref{sec:dispmanage} we assume an optical communication line of 
$100\,\mathrm{km}$ length with dispersion parameters $|\beta_2^{(1,2)}|=0.5\,\mathrm{ps^2 km^{-1}}$ and $|\beta_3^{(1,2)}|=0.07\,\mathrm{ps^3 km^{-1}}$, 
which all change sign every $2\,\mathrm{km}$, and $\alpha_{(1,2)}=0$, $\gamma_{(1,2}=0.1\,\mathrm{W/m}$, and $T_R=3\,\mathrm{fs}$. As before, the 
parameters of the two unchirped Gaussian pulses are $P_0^{(1)}=0.625\,\mathrm{mW}$, $P_0^{(2)}=0.3125\,\mathrm{mW}$, and $T_0^{(1,2)}=80\,\mathrm{fs}$.
The wavelengths are again $\lambda_0^{(1)}=1550\,\mathrm{nm}$ and $\lambda_0^{(1)}=1300\,\mathrm{nm}$. The group velocity mismatch
is $\delta = 0.015625\,\mathrm{fs/m}$ and cross-phase modulation parameters are $B_{1,2}=C_{1,2}=2$. The same computational
parameters are used as in Sect.~\ref{sec:dispmanage}: The temporal window has the width  $[-30\,\mathrm{ps},30\,\mathrm{ps}-\Delta T]$
and $N=4096$, $h=40\,\mathrm{m}$ are applied. 

\begin{figure}[t]
\centering\includegraphics[angle=-90,width=6cm]{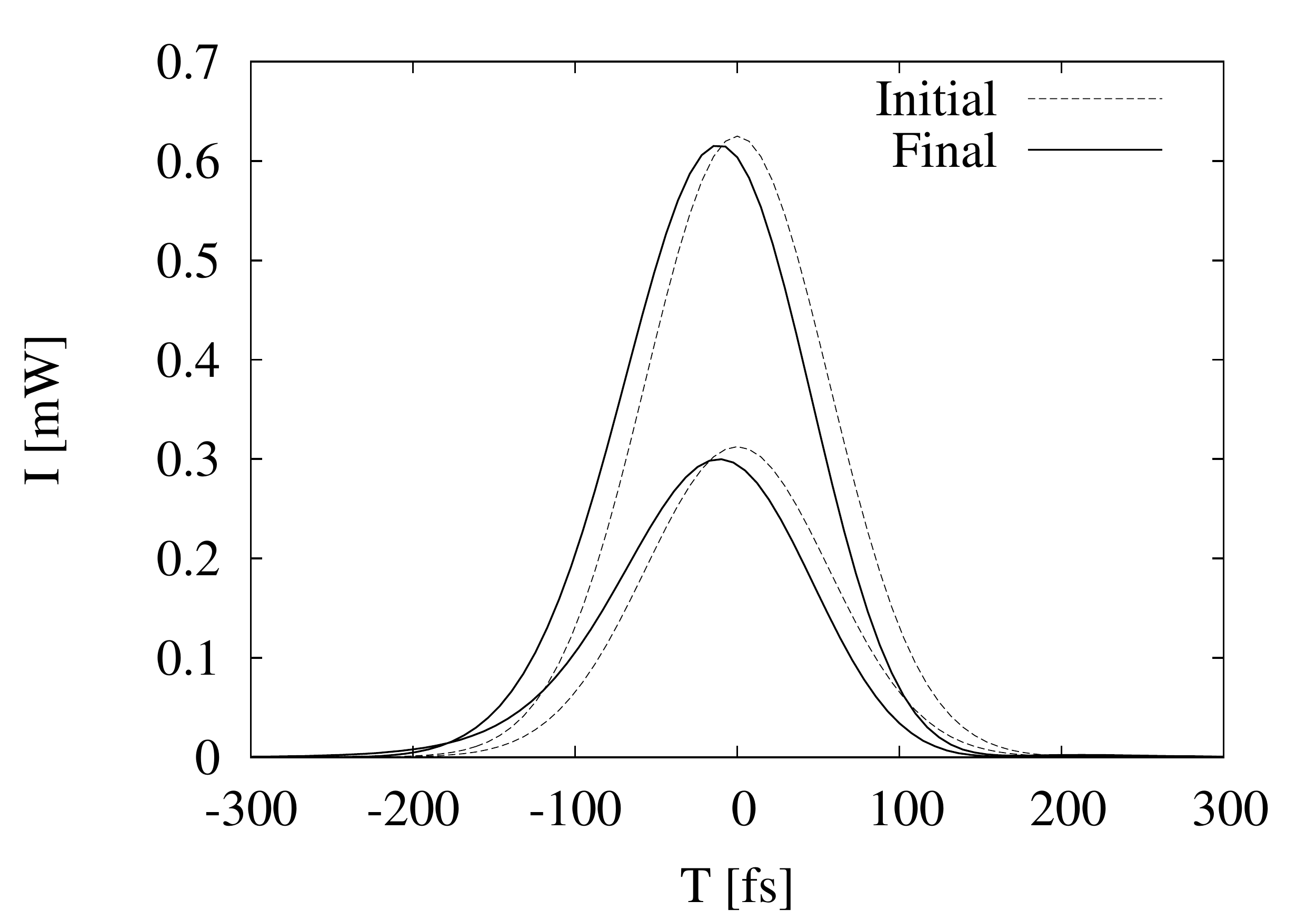}\includegraphics[angle=-90,width=6cm]{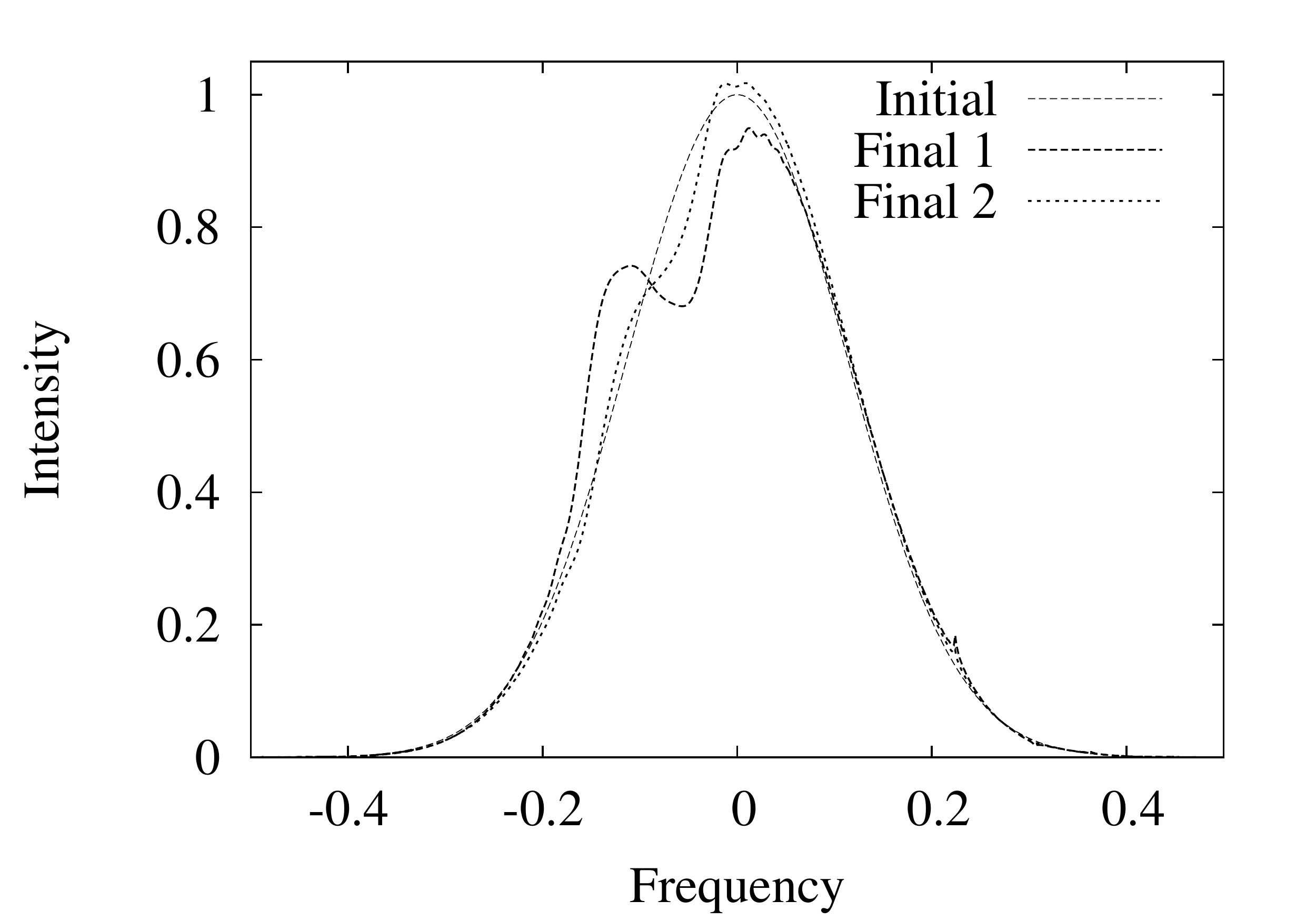}

\caption{Pulse shape and spectrum of two coupled pulses after propagating $100\,\mathrm{km}$ through the idealized dispersion-managed fiber of
Benchmark 2.} 
\label{pic:bench4} 
\end{figure}

During propagation both pulses are experiencing almost undisturbed soliton-like oscillations every $4\,\mathrm{km}$.  
Figure~\ref{pic:bench4} compares the final signal shapes and spectra with the respective initial ones, where Pulse 2 has been shifted for 
visualization by $-1.5625\,\mathrm{ps}$. Both pulses are delayed by roughly $10\,\mathrm{fs}$ but the signal shape is quite well preserved;
the spectral alteration being rather moderate in both cases. In the left graphic of Fig.~\ref{pic:bench4} Pulse 1 and 2 are easily distinguished;
in the right graphic the final spectra of Pulse 1 and 2 are specially indicated. 

Finally, we comment on typical run times of the proposed split-step Fourier methods. Our implementation is in FORTRAN 90 and uses the Netlib 
NAPACK Fast Fourier Transformation (FFT) routines, which are coded in FORTRAN 77 \cite{Hager-88}. Compiled with usual optimizations, 
the two-mode computation of Fig.~\ref{pic:bench3} required $\sim 48$ seconds on a single Intel Xeon E5 CPU with $2.1\,\mathrm{GHz}$. Dependence
on the number of Fourier modes $N$ as well as the number of spatial steps $M$ is linear and each computation of the convergence analysis of 
Fig.~\ref{pic:bench3b} is therefore four times more expensive than the next coarser one. On the same CPU, the dispersion-managed two-mode simulation of 
Fig.~\ref{pic:bench4} ran for $\sim 100$ seconds, its single-mode analogue of Fig.~\ref{pic:bench2b} required $\sim 40$ seconds. These moderate
run times and the given results provide evidence for the relevance of the proposed numerical methods for practical long-distance fiber optical 
communication line design.

\section{Conclusions}\label{sec:conclusions}
Reliable extensions of the classical SSFM into the regime of ultra-fast pulses have been derived and demonstrated
for typical Gaussian communication pulses in highly nonlinear and dispersion-managed long-distance optical fibers. The primary difficulty in this regime
lies in the appropriate mathematical treatment of the additional nonlinear terms modeling signal self-steepening and stimulated Raman scattering. For
the case of the single-mode equation (\ref{eq:pulsefast}) and the two-mode system (\ref{eq:tmodeul}) it was shown that under Madelung
transformation all nonlinearities can be effectively combined into an inhomogeneous system of advection equations of the signal intensities 
and phases. Following upwind and slope-limiting ideas, originally developed in the context of supersonic hydrodynamics, a robust numerical 
method is then derived for the single-mode nonlinear sub-operator and incorporated into a symmetric SSFM. Reliable convergence and
numerical approximation accuracy of second order is demonstrated for the overall method. While it would be principally feasible to apply the exact
same approach to the two-mode case and the correspondingly derived four-dimensional system (\ref{eq:tmhyp}), we have opted for now for a 
mathematically less involved fractional step approach and apply two single-field nonlinear sub-operators successively to approximate the solution
of (\ref{eq:tmhyp}). This single-field sub-operator is derived as a straightforward extension of the slope-limited upwind method for 
the single-mode case. Incorporated into a two-mode SSFM, the overall numerical scheme converges reliably, yet, in a highly nonlinear test case 
only an order of accuracy of $\sim 1.5$ is obtained. Future work will concentrate on developing an unsplit scheme for (\ref{eq:tmhyp}). 
It is expected that such a method should obtain an order of accuracy close to 2 while being of comparable computational expense and robustness as 
the two-mode SSFM proposed in here.

\begin{acknowledgement}
This work was supported by the Department of Defense and used resources of the Extreme Scale
Systems Center at Oak Ridge National Laboratory.
\end{acknowledgement}

\bibliographystyle{spmpsci}
\bibliography{../../Literature/all}

\end{document}